\documentclass[12pt]{amsart}
\usepackage{amscd,amsmath,amssymb,amsfonts}
\usepackage[cmtip, all]{xy}

\newtheorem{thm}{Theorem}
\newtheorem{prop}[thm]{Proposition}
\newtheorem{lem}[thm]{Lemma}
\newtheorem{cor}[thm]{Corollary}

\theoremstyle{definition}

\theoremstyle{remark}

\numberwithin{equation}{section}

\newcommand{\sE}{{\mathcal E}}

\newcommand{\sH}{{\mathcal H}}
\newcommand{\sI}{{\mathcal I}}

\newcommand{\sM}{{\mathcal M}}

\newcommand{\sO}{{\mathcal O}}

\newcommand{\sV}{{\mathcal V}}

\newcommand{\A}{{\mathbb A}}

\newcommand{\C}{{\mathbb C}}

\renewcommand{\P}{{\mathbb P}}
\newcommand{\Q}{{\mathbb Q}}

\newcommand{\Z}{{\mathbb Z}}

\renewcommand{\L}{{\mathbb{L}}}

\newcommand{\CH}{{\rm CH}}

\newcommand{\Pic}{{\rm Pic}}
\newcommand{\Div}{{\rm Div}}

\newcommand{\Spec}{{\rm Spec \,}}

\newcommand{\0}{\emptyset}

\newcommand{\id}{{\operatorname{Id}}}

\newcommand{\Sch}{{\operatorname{\mathbf{Sch}}}}

\newcommand{\Sm}{{\mathbf{Sm}}}

\newcommand{\GL}{{\operatorname{\rm GL}}}

 \newcommand{\Ab}{{\mathbf{Ab}}}

\newcommand{\Sym}{{\operatorname{\rm Sym}}}

\newcommand{\ds}{{/\kern-3pt/}}

\newcommand{\lci}{\text{l.c.i.}}

\newcommand{\Tor}{{\operatorname{Tor}}}
\newcommand{\Proj}{{\operatorname{Proj}}}

\newcommand{\pt}{{\operatorname{pt}}}
\newcommand{\oo}{\otimes}
\newcommand{\cn}{{\tilde{c}_1}}

\newcommand{\ZZ}{{\mathsf{Z}}}

\author{M. Levine and R. Pandharipande}

\date{ 3 May 2006}

\begin{document}


\title{Algebraic cobordism revisited}

\begin{abstract}
We define a cobordism theory in algebraic geometry
based on normal crossing degenerations with 
double point singularities. 
The main result is the equivalence of double point cobordism
 to the theory of algebraic cobordism previously defined by Levine and
 Morel. Double
 point cobordism provides a simple, geometric presentation
 of algebraic cobordism theory.
 As a corollary,
 the Lazard ring given by products of projective spaces
 rationally generates all nonsingular projective varieties
 modulo double point degenerations.
 Double point  degenerations 
  arise naturally in 
  relative Donaldson-Thomas theory. We use double point
  cobordism to prove all the degree 0 conjectures
  in Donaldson-Thomas theory: absolute, relative, and
  equivariant.\end{abstract}

 \maketitle

 
\setcounter{section}{-1} 
\section{Introduction}

\subsection{Overview}
A first idea for defining cobordism in algebraic geometry is to impose the
relation
\begin{equation}\label{naive}
[\pi^{-1}(0)] = [\pi^{-1}(\infty)]
\end{equation}
for smooth fibers of a projective morphism
\begin{equation*}
\pi:Y \to \P^1.
\end{equation*}
The resulting theory bears no resemblance to complex cobordism.

A successful theory of algebraic cobordism 
has been constructed in \cite{ACI,ACII}
from Quillen's axiomatic perspective. The goal is to define a
universal oriented Borel-Moore cohomology theory of schemes. 
An introduction to 
algebraic cobordism can be found in \cite{ICM1, SAC, CR1,CR2}.

A second idea for defining algebraic cobordism geometrically is to
impose relations obtained by fibers
of $\pi$ with  
 normal crossing singularities.
The simplest of these are the {\em double point degenerations} ---
where the fiber is a union of two smooth transverse divisors. 
We prove the cobordism theory obtained from
double point degenerations {\em is} algebraic cobordism.

Algebraic cobordism may thus be viewed both functorially and geometrically.
In practice, the different perspectives are very useful.
We prove several conjectural formulas concerning the virtual class of the
Hilbert scheme of points of a 3-fold as an application.

\subsection{Schemes and morphisms} 
Let $k$ be a field of characteristic 0.
Let $\Sch_k$ be the category of 
 separated schemes of finite type over $k$, and
let $\Sm_k$ be the full subcategory of smooth quasi-projective schemes.

For $X\in \Sch_k$,
let $\sM(X)$ denote the set of isomorphism classes over $X$ of
projective morphisms
\begin{equation} \label{fyx}
f:Y \rightarrow X
\end{equation}
with $Y\in \Sm_k$. The set $\sM(X)$ is a monoid under
disjoint union of domains. Let $\sM(X)^+$ denote
the group completion of $\sM(X)$. 

Alternatively, $\sM(X)^+$ is the
free abelian group generated by morphisms \eqref{fyx}
where $Y$ is irreducible. Let 
$$[f:Y \to X] \in \sM(X)^+$$
denote the element determined by the morphism.

\subsection{Double point degenerations}
Let $Y\in \Sm_k$ be of pure dimension. A morphism
$$\pi:Y\to \P^1$$
is a {\em double point degeneration}
over  $0\in \P^1$ if
$$\pi^{-1}(0)=A \cup B $$ 
where 
$A$ and $B$ are smooth Cartier divisors intersecting transversely in $Y$. 
The intersection $$D=A\cap B$$ is the {\em double point locus} of $\pi$
over $0\in \P^1$.

Let $N_{A/D}$ and $N_{B/D}$ denote the normal bundles of $D$ in $A$ and $B$
respectively.
Since 
$O_D(A+B)$ is trivial,
$$N_{A/D} \otimes N_{B/D} \cong O_D.$$
Hence, the projective bundles
\begin{equation}\label{pbd}
\P(O_D\oplus N_{A/D})\to D \ \ \text{and} \ \ \P(O_D\oplus N_{B/D})\to D
\end{equation}
are isomorphic. Let 
$$\P(\pi) \rightarrow D$$
denote either of \eqref{pbd}.

\subsection{Double point relations}
Let $X\in \Sch_k$, and let 
$p_1$ and $p_2$ denote the projections to the first and second factors of 
 $X \times \P^1$ respectively.

Let $Y\in \Sm_k$ be of pure dimension. Let
$$\pi: Y \to X \times \P^1$$ be a projective morphism
for which the 
composition
\begin{equation}\label{pi2}
\pi_2 = p_2\circ \pi: Y \to \P^1
\end{equation}
is 
a double point degeneration over $0 \in \P^1$. 
Let 
$$[A\to X],\ [B\to X], \  [\P(\pi_2) \to X] \in \sM(X)^+$$
be obtained from the fiber $\pi_2^{-1}(0)$ and the
morphism $p_1\circ \pi$.

For each
regular value 
$\zeta\in \P^1(k)$ of $\pi_2$, define an associated
{\em double point relation} over $X$ by
\begin{equation}\label{dpr}
[Y_\zeta \to X] -[A \to X] - [B \to X] + [\P(\pi_2) \to X] 
\end{equation}
where $Y_\zeta = \pi_2^{-1}(\zeta)$.

Let $\mathcal{R}(X)\subset \sM(X)^+$ be the subgroup generated
by {\em all} double point relations over $X$.

\subsection{Naive cobordism}
Naive cobordism \eqref{naive} may be viewed as a special case
of a double point relation.

Let $Y\in \Sm_k$ be of pure dimension.
Let
$$\pi: Y \to X \times \P^1$$ be a
projective morphism with $\pi_2=p_2\circ \pi$ 
smooth over $0,\infty\in \P^1$.
We may view $\pi_2$ as a double point degeneration
over $0\in \P^1$ with
$$\pi_2^{-1}(0)= A \cup \emptyset.$$
The associated double point relation is
$$[Y_\infty \to X] - [Y_0 \to X]\in \mathcal{R}(X).$$

\subsection{Algebraic cobordism} 
\label{dpc}
The central object of the paper is the quotient
$$\omega_*(X) = \sM(X)^+/\mathcal{R}(X)$$ 
defining a double point cobordism theory.
Let $\Omega_*(X)$ be the theory of algebraic cobordism
defined in \cite{ACI, ACII}. 

\begin{thm} \label{one}
There is a canonical isomorphism $\omega_*(X) \cong \Omega_*(X)$.
\end{thm}

Theorem \ref{one} may be viewed as a geometric
presentation of $\Omega_*(X)$ via the simplest possible 
cobordisms. 
A homomorphism 
\begin{equation}\label{hom}
\omega_*(X) \to \Omega_*(X)
\end{equation} is obtained
immediately from the definitions once the double point
relations are shown to hold in $\Omega_*(X)$.
The inverse is more difficult
to construct.

\begin{thm} \label{two} 
$\omega_*$ determines an oriented Borel-Moore functor of geometric type on
$\Sch_k$.
\end{thm}

Since algebraic cobordism is the universal Borel-Moore
functor of geometric type on $\Sch_k$, an
inverse $$\Omega_*(X)\to \omega_*(X)$$
to \eqref{hom} is 
obtained from Theorem \ref{two}.

Oriented theories and Borel-Moore functors are discussed
in Sections \ref{ot}- \ref{gggg} following \cite{ACI,ACII}.
The proof of Theorem \ref{two}, presented in Sections \ref{fff}-\ref{eee},
is the technical heart
of the paper. 
The key geometric step is the construction of a
formal group law for $\omega_*$ in Section \ref{rrr}.
Theorem \ref{one} is proven in Section \ref{t1}.

\subsection{Algebraic cobordism over a point} 
Denote $\text{Spec}(k)$ by $k$. 
Let $\L_*$ be the Lazard ring \cite{Lazard}.
The canonical map
$$\L_* \to \Omega_*(k)$$
is proven to be an isomorphism in \cite{ACI,ACII}.
By Theorem \ref{one},
$$\L_* \cong \omega_*(k).$$
A basis of $\omega_*(k)\otimes_{\mathbb Z} {\mathbb Q}$
is formed by the products of projective spaces.

\begin{cor} 
\label{bas}
We have 
$$\omega_*(k)\otimes_{\mathbb Z} {\mathbb Q} 
=  \bigoplus_{\lambda} \  {\mathbb Q} [\P^{\lambda_1} \times ... 
\times \P^{\lambda_{\ell(\lambda)}}]$$
where the sum is over all partitions $\lambda$.
\end{cor}

\subsection{Donaldson-Thomas theory}\label{dt}
Corollary \ref{bas} is directly applicable 
to the Donaldson-Thomas theory of 3-folds.

Let $X$ be a smooth projective 3-fold over $\mathbb{C}$, and
let $\text{Hilb}(X,n)$ be the Hilbert scheme of $n$ points.
Viewing the Hilbert scheme as the moduli space of ideal sheaves
$I_0(X,n)$, a natural
0-dimensional virtual Chow class can be constructed
$$[\text{Hilb}(X,n)]^{vir} \in A_0(\text{Hilb}(X,n),\mathbb{Z}),$$
see \cite{MNOP1,MNOP2,T}. The degree 0 Donaldson-Thomas invariants
are defined by 
$$N^X_{n,0} = \int_{[\text{Hilb}(X,n)]^{vir}} 1.$$
Let  
$$\ZZ(X,q) = 1 + \sum_{n\geq 1} N^X_{n,0}\ q^n$$
be the associated partition function.

\vspace{10pt}

\noindent{\bf Conjecture 1.} \cite{MNOP1}\,  $\ZZ(X,q) =
M(-q)^{\int_X c_3(T_X \otimes K_X)}.$

\vspace{10pt}

Here, $M(q)$ denotes the MacMahon function,
$$M(q) = \prod_{n\geq 1} \frac{1}{(1-q^n)^n},$$
the generating function of 3-dimensional partitions  \cite{Stanley}.

For a nonsingular divisor $S\subset X$, a relative
Donaldson-Thomas theory\footnote{See \cite{MNOP1,OP}
for a discussion. A full foundational treatment of the relative
theory has not yet appeared.}
is defined via the moduli space of
relative ideal sheaves $I_0(X/S,n)$.
The degree 0 relative invariants, 
$$N^{X/S}_{n,0} = \int_{[I_0(X/S,n)]^{vir}} 1,$$
determine a relative partition function
$$\ZZ(X/S,q) = 1 + \sum_{n\geq 1} N^{X/S}_{n,0}\ q^n.$$

Let $\Omega_X[S]$ denote the locally free sheaf of
differential forms of $X$ with logarithmic poles
along $S$. Let 
$$T_X[-S]= \Omega_X[S]\ ^\vee,$$
denote the dual sheaf of tangent fields with logarithmic zeros.
Let 
$$K_X[S]= \Lambda^3 \Omega_X[S]$$
denote the logarithmic canonical class.

\vspace{10pt}

\noindent{\bf Conjecture 2.} \cite{MNOP2}\,  $\ZZ(X/S,q) =
M(-q)^{\int_X  c_3(T_X[-S] \otimes K_X[S])}.$

\vspace{10pt}

We prove Conjectures 1 and 2. 
An equivariant version of Conjecture 1 proposed in \cite{bp} is
also proven.
Corollary \ref{bas} 
reduces
the results to toric cases previously calculated in
\cite{MNOP1,MNOP2}. 
The proofs are presented in Section \ref{last}.

\subsection{Double point relations in DT theory}
Double point relations
naturally arise in degree 0 Donaldson-Thomas theory by the
following construction.

Let $Y\in \Sm_{\mathbb{C}}$ be a 4-dimensional projective
variety, and let
$$\pi:Y \to \P^1$$ be a double point degeneration 
over  $0\in \P^1$. Let 
$$\pi^{-1}(0)=A\cup B.$$
The degeneration formula in relative Donaldson-Thomas
theory yields
\begin{equation}\label{x}
\ZZ(Y_\zeta) = \ZZ(A/D)\cdot  \ZZ(B/D)
\end{equation}
for a $\pi$-regular value $\zeta\in \P^1$, see \cite{MNOP2}.

Since the deformation to the normal cone of $D \subset A$ is 
a double point degeneration,
\begin{equation}\label{xx}
\ZZ(A) = \ZZ(A/D)\cdot  \ZZ(\P(O_D \oplus N_{A/D})/D).
\end{equation}
On the right, the divisor $D\subset \P(O_D \oplus N_{A/D})$
is included with normal bundle $N_{A/D}$.
Similarly,
\begin{equation}\label{xxx}
\ZZ(B) = \ZZ(B/D)\cdot  \ZZ(\P(O_D \oplus N_{B/D})/D)
\end{equation}
where the divisor 
$D\subset \P(O_D \oplus N_{A/D})$ is included with
normal bundle $N_{B/D}$.

Since $N_{A/B} \otimes N_{B/D} \cong O_D$, the deformation
of $\P(O_D \oplus N_{A/D})$ to the normal cone of
$D\subset \P(O_D \oplus N_{A/D})$ yields
$$
\ZZ(\P(\pi)) = 
\ZZ(\P(O_D \oplus N_{A/D})/D)\cdot  
\ZZ(\P(O_D \oplus N_{B/D})/D).
$$
When combined with equations \eqref{x}-\eqref{xxx}, we find
\begin{equation}\label{dpdt}
\ZZ(Y_\zeta) \cdot \ZZ(A)^{-1}\cdot \ZZ(B)^{-1}\cdot \ZZ(\P(\pi)) =1
\end{equation}
which is the double point relation \eqref{dpr} over  $\text{Spec}
(\mathbb{C})$
in multiplicative
form.

\subsection{Gromov-Witten speculations}
Let $X$ be a nonsingular projective variety over $\mathbb{C}$.
Gromov-Witten theory concerns integration against
the virtual class,
$$[\overline{M}_{g,n}(X,\beta)]^{vir} \in 
H_*(\overline{M}_{g,n}(X,\beta), \mathbb{Q}),$$
of the moduli space of 
stable maps to $X$.
 
There are two main techniques available in Gromov-Witten theory:
localization \cite{GP,Kon} and degeneration \cite{EGH,IP,LR,L,mp}. 
Localization is most effective
for toric targets --- all the Gromov-Witten data
of products of projective spaces are accessible by localization.
The degeneration formula yields Gromov-Witten relations precisely
for double point degenerations.

By Corollary \ref{bas}, {\em all} varieties are
linked to products of projective spaces by double point
degenerations. We can expect, therefore, that
many aspects of the Gromov-Witten theory of arbitrary varieties 
will follow the behavior found in toric targets.
An example is the following speculation about the virtual
class --- which, at present,
appears  out of reach of Corollary \ref{bas}.

\vspace{10pt}

\noindent{\bf Speculation.}  {\em The push 
forward $\epsilon_*[\overline{M}_{g,n}(X,\beta)]^{vir}$
 via the canonical map
$$\epsilon: \overline{M}_{g,n}(X,\beta) \to \overline{M}_{g,n}$$
lies in the tautological ring 
$$RH^*(\overline{M}_{g,n},\mathbb{Q}) 
\subset H^*(\overline{M}_{g,n},\mathbb{Q}).$$}
\vspace{5pt}

See \cite{FPRT, pan} for a discussion of similar (and stronger)
statements. In particular, a definition of the
tautological ring can be found there.

Gromov-Witten theory is most naturally viewed as an aspect
of symplectic geometry. The construction of
a parallel symplectic cobordism theory based on double
point degenerations appears to be a natural path to follow.

\subsection{Acknowledgments}
We thank D. Maulik, A. Okounkov, and B. Totaro for useful discussions about
double point degenerations, Gromov-Witten theory, and
algebraic cobordism.

Conjecture 1 has been recently proven by J. Li \cite{lic1}.
Li's method is to show $\ZZ(X,q)$ depends only upon the
Chern numbers of $X$
by an explicit (topological) study of the cones
defining the virtual class. The result is then obtained
from the toric calculations of \cite{MNOP2}
via the complex cobordism class.
A proof of Conjecture 1 in case $X$ is a Calabi-Yau 3-fold via a
study of self-dual obstruction theories appears in
\cite{B,BF}. Our proof is direct and algebraic,
but depends upon the construction of
relative Donaldson-Thomas theory (which is required
in any case for the calculations of \cite{MNOP2}).

M. L. was 
supported by the Humboldt Foundation through the 
Wolfgang Paul Program  and the NSF 
via grants DMS-0140445 and  DMS-0457195.
R. P. was supported by the Packard foundation and
the NSF via grant DMS-0500187. 
The research was partially pursued during a visit
of R. P. to the Instituto Superior T\'ecnico in Lisbon
in the fall of 2005.

\section{Oriented theories} \label{ot}
\subsection{$\Omega_*$} 
Theorem \ref{one} is proven 
for algebraic cobordism $\Omega_*$ viewed
as an oriented Borel-Moore
homology theory on $\Sch_k$. 
We start by reviewing the  
definitions of oriented homology and cohomology theories following
\cite{ACI,ACII}.

\subsection{Notation}
Let $X\in \Sch_k$. A {\em divisor} $D$ on $X$ will be understood
to be Cartier unless otherwise stated. The line bundle associated
to the locally free
sheaf $\sO_X(D)$ is denoted $O_X(D)$.

Let $\sE$ be a rank $n$ locally free sheaf $\sE$ on $X$. Let 
$$q:\P(\sE)\to X$$ denote the 
projective bundle
$\Proj_X(\Sym^*(\sE))$ of rank one 
{\em quotients} of $\sE$ with tautological quotient 
invertible sheaf 
$$q^*\sE\to\sO(1)_\sE.$$ We let $O(1)_\sE$ denote the 
line bundle on $\P(\sE)$ with sheaf of 
sections $\sO(1)_\sE$. 
The subscript ${}_\sE$ is omitted
if the context makes the meaning clear. The notation $\P_X(\sE)$ 
is used to emphasize the base scheme $X$.

Two morphisms $f:X\to Z$, $g:Y\to Z$ in $\Sch_k$ are {\em Tor-independent} if, 
for each triple of points $$x\in X,\ \ y\in Y,\ \ z\in Z$$ 
satisfying $f(x)=g(y)=z$, 
\[
\Tor_p^{\sO_{Z,z}}(\sO_{X,x},\sO_{Y,y})=0
\]
for $p>0$.

A closed immersion $i:Z\to X$ in $\Sch_k$ is a {\em regular embedding} 
if the ideal sheaf $\sI_Z$ is locally generated by a regular sequence.  
A morphism $f:Z\to X$ in $\Sch_k$ is {\em l.c.i.} if 
$$f=p\circ i$$ 
where $i:Z\to Y$ is a regular embedding and $p:Y\to X$ is a smooth  
morphism.\footnote{For us, a smooth 
morphism is smooth {\em{and}} quasi-projective.}
 L.c.i. morphisms are closed under composition. 

If $f:Z\to X$ and $g:Y\to X$ are 
Tor-independent  morphisms in $\Sch_k$ and $f$ is an 
l.c.i.-morphism, then 
$$p_1:Z\times_XY\to Y$$ is an l.c.i. morphism.

For a full subcategory $\sV$ of $\Sch_k$, 
let $\sV'$ denote the category with 
$$\text{Ob}(\sV')=\text{Ob}(\sV)$$
and arrows given by {projective} morphisms of
schemes. 

Let $\Ab_*$ denote the category of graded abelian groups.
A functor 
$$F:\Sch_k'\to \Ab_*$$ is {\em additive} if $F(\0)=0$ and the canonical map 
$$F(X)\oplus F(Y)\to F(X\coprod Y)$$ is an 
isomorphism for all $X$, $Y$ in $\Sch'_k$.

\subsection{Homology}
We review 
the definition of an oriented Borel-Moore homology theory from \cite{ACII}.
We refer the reader to \cite{ACII} for a more leisurely discussion.

An {\em oriented 
Borel-Moore homology theory} $A_*$ on
$\Sch_k$ consists of the following data:
\begin{enumerate}
\item[(D1)] An additive functor
\[
A_* : \Sch_k' \to \Ab_*\, , \, X\mapsto A_*(X).
\]
\item[(D2)] For each \lci\ morphism $f:Y\to X$ 
in $\Sch_k$ of relative dimension $d$, a homomorphism of graded groups 
\[
f^*:A_*(X)\to A_{*+d}(Y).
\]
\item[(D3)] 
For each pair
$(X,Y)$  in $\Sch_k$,
a bilinear graded pairing
\begin{align*}
A_*(X)\oo A_*(Y) &\to A_*(X\times_k Y)\\
u\oo v&\mapsto u\times v,
\end{align*}
which 
is commutative, associative,  and admits a distinguished
element $1\in A_0(\text{Spec}(k))$ as a unit.
\end{enumerate}
The pairing in (D3) is the {\em external product}.
The data (D1)-(D3) are required to satisfy six conditions:
\begin{enumerate}
\item[(BM1)] 
Let  
$f:Y\to X$ and $g:Z\to Y$ be \lci\ morphisms
in $\Sch_k$ of pure relative dimension. Then,
$$(f\circ g)^*=g^*\circ f^*.$$
Moreover, $\id_X^*=\id_{A_*(X)}$.
\item[(BM2)] Let $f:X\to Z$ and $g:Y\to Z$ be 
$\Tor$-independent morphisms in $\Sch_k$ where
 $f$ is projective and $g$ is \lci\  In the  
 cartesian square
\[
\xymatrix{
W\ar[r]^{g'}\ar[d]_{f'}&X\ar[d]^f\\
Y\ar[r]_g&Z\, ,}
\]
$f'$ is projective and $g'$ is  \lci\  Then, 
$$g^*f_*=f'_*g^{\prime*}.$$
\item[(BM3)] Let $f:X'\to X$  and $g:Y'\to Y$ be morphisms in 
$\Sch_k$.\\ If $f$ and $g$ are projective, 
then  
\[
(f\times g)_* (u'\times v') = f_*(u')\times g_*(v').
\]
for $u'\in A_*(X')$ and $v'\in A_*(Y')$.\\ 
If $f$ and $g$ are \lci, then 
\[
(f\times g)^* (u\times v) = f^*(u)\times g^*(u')
\]
for $u\in A_*(X)$ and 
$v\in A_*(Y)$.
\item[(PB)] For a line bundle $L$ on $Y\in\Sch_k$ with zero 
section $$s:Y\to
L,$$ define the operator
\[
\cn(L):A_*(Y)\to A_{*-1}(Y)
\]
by $\cn(L)(\eta)=s^*(s_*(\eta))$.  \\ \\
Let $\sE$ be a rank $n+1$ locally free 
sheaf on $X\in\Sch_k$, with associated
projective bundle $$q:\P(\sE)\to X.$$ For 
$i=0,\dots,n$, let 
\[
\xi^{(i)}: A_{*+i-n}(X)\to A_{*}(\P(\sE))
\]
be the composition of $$q^*: A_{*+i-n}(X)\to A_{*+i}(\P(\sE))$$ followed by 
$$\cn(O(1)_\sE)^i: A_{*+i}(\P(\sE))\to A_{*}(\P(\sE)).$$ 
Then the homomorphism
\[
\Sigma_{i=0}^{n-1} \xi^{(i)}: \oplus_{i=0}^n A_{*+i-n}(X)\to A_*(\P(\sE))
\]
is an isomorphism.
\item[(EH)] Let $E\to X$ be a vector bundle of rank
$r$ over $X\in\Sch_k$, 
and let $p:V\to X$ be an $E$-torsor. Then 
$$p^*:A_*(X)\to A_{*+r}(V)$$ is an 
isomorphism.
\item[(CD)] For  integers $r,N>0$, let 
$$W=\underbrace{\P^N\times_S\ldots\times_S\P^N}_r,$$ 
and let
$p_i:W\to\P^N$ be the $i$th projection. Let $X_0,\ldots, X_N$ be the
standard homogeneous coordinations on $\P^N$, let
$n_1,\ldots, n_r$ be non-negative 
integers, and let $i:E\to W$ be the subscheme defined
by
$\prod_{i=1}^rp_i^*(X_N)^{n_i}=0$.  Then $$i_*:A_*(E)\to
A_*(W)$$ is injective.
\end{enumerate}

Comments about (CD) in relation to a more natural filtration
condition can be found in \cite{ACII}

The most basic example of an oriented Borel-Moore homology theory on
$\Sch_k$ is the 
 Chow group functor 
$$X\mapsto \CH_*(X)$$ with  
projective push-forward
and l.c.i.
pull-back 
given by Fulton \cite{Fulton}.

\label{axhom}

\subsection{Cohomology}
Oriented  cohomology theories on $\Sm_k$
are defined axiomatically in \cite{ACI}.
The axioms are very similar to those discussed in
Section \ref{axhom}.

An oriented cohomology theory $A^*$ on $\Sm_k$ 
can be obtained from 
an oriented  Borel-Moore homology theory $A_*$ on $\Sch_k$ by reindexing.
If $X\in \Sm_k$ is irreducible,
\[
A^*(X)=A_{\dim X-*}(X).
\]
In the reducible case, the reindexing is applied to each component
via the additive property.

$A^*(X)$ is a commutative graded ring with unit. The product is defined by 
$$a\cup b=\delta^*(a\times b)$$
where $\delta:X\to X\times X$ is the diagonal. The unit 
is 
$$1_X=p_X^*(1)$$ 
where $p_X:X\to\text{Spec}(k)$ is the structure morphism.

The first Chern class has the following interpretation in oriented
cohomology. Let $L$ be a line bundle on $X$, and 
let $$c_1(L)=\cn(L)(1_X)\in A^1(X),$$ then
\[
\cn(L)(a)=c_1(L)\cup a
\]
for all $a\in A^*(X)$.

Let $f:Y\to X$ be a morphism in $\Sch_k$  with $X\in\Sm_k$. 
Then 
$$(f,\id):Y\to X\times Y$$ is a regular embedding. The pairing
\begin{align*}
A^m(X)\otimes A_n(Y)&\to A_{n-m}(Y)\\
a\otimes b&\mapsto (f,\id)^*(a\times b)
\end{align*}
makes $A_*(Y)$ a graded $A^*(X)$-module (with $A_{-n}(Y)$ in degree $n$).

\section{Algebraic cobordism theory $\Omega_*$}
\subsection{Construction}
Algebraic cobordism theory
is constructed in \cite{ACI}, and many
 fundamental properties of $\Omega_*$ are
verified there.  
The program is completed in \cite{ACII} by proving
$\Omega_*$ 
is a universal
oriented Borel-Moore homology theory on $\Sch_k$.
The result requires the construction of 
pull-back maps for l.c.i. morphisms.
We give a basic sketch of $\Omega_*$ here.

\subsection{$\underline{\Omega}_*$}
For $X\in\Sch_k$, $\Omega_n(X)$ is generated 
(as an abelian group) by {\em cobordism cycles} 
$$(f:Y\to X, L_1,\ldots, L_r),$$
where $f$ is a projective morphism, $Y\in\Sm_k$ is irreducible of 
dimension $n+r$ over $k$, 
and the $L_i$ are line bundles on $Y$. 
We identify two cobordism cycles if they are isomorphic over $X$ 
up to reorderings of the line bundles $L_i$. 

We will impose several relations on cobordism cycles.
To start, two basic relations are imposed:
\begin{enumerate}
\item[I.] 
If there exists a smooth morphism $\pi:Y\to Z$ and 
 line bundles $M_1,\ldots, M_{s>\dim_k Z}$ on $Z$ with 
$L_i\cong \pi^*M_i$ for  $i=1,\ldots, s\le r$, then
$$(f:Y\to X, L_1,\ldots, L_r)=0.$$ 
\item[II.] If $s:Y\to L$ is a section of a line bundle
 with smooth associated divisor 
$i:D\to Y$, then  
$$\ \ \ \ \ \ \ \ 
(f:Y\to X, L_1,\ldots, L_r, L)=(f\circ i:D\to X, i^*L_1,\ldots, i^*L_r).$$
\end{enumerate}
The graded group generated by cobordism cycles modulo relations
I and II is
denoted
$\underline{\Omega}_*(X)$.

Relation II yields as a special case the naive cobordism relation. 
Let $$\pi:Y\to X\times\P^1$$ be a projective morphism with $Y\in\Sm_k$
for which 
$p_2\circ \pi$ is  transverse to the inclusion $\{0,\infty\}\to \P^1$. 
Let $L_1,\ldots, L_r$ be line bundles on $Y$, and let 
$$i_0:Y_0\to Y, \ \ i_\infty:Y_\infty\to Y$$ be 
the inclusions of the fibers over $0,\infty$. Then 
\[
(p_1\circ \pi:Y_0\to X, i_0^*L_1,\ldots, i_0^*L_r)
=(p_1\circ \pi :Y_\infty\to X, 
i_\infty^*L_1,\ldots, i_\infty^*L_r)
\]
in $\underline{\Omega}_*(X)$.

Several structures are easily constructed on $\underline{\Omega}_*$.
For a projective morphism $g:X\to X'$, define
$$g_*:\underline{\Omega}_*(X) \to \underline{\Omega}_*(X')$$
by the rule 
\[
g_*(f:Y\to X, L_1,\ldots, L_r)=(g\circ f:Y\to X', L_1,\ldots, L_r).
\]
Similarly 
evident pull-backs for smooth morphisms and external products
exist for $\underline{\Omega}_*$.

The 
Chern class operator $\cn(L):\underline{\Omega}_n(X)\to 
\underline{\Omega}_{n-1}(X)$ is defined by the following formula:
\[
\cn(L)((f:Y\to X, L_1,\ldots, L_r))=(f:Y\to X, L_1,\ldots, L_r,f^*L).
\]

\subsection{$\Omega_*$}
 Contrary to the purely 
topological theory of complex cobordism, relations I
 and II do not suffice to define $\Omega_*$. 
One needs to impose the  formal group law.

A (commutative, rank one)  {\em formal group law} 
over a commutative ring $R$ is a power 
series $F(u,v)\in R[[u,v]]$ satisfying the formal 
relations of identity, commutativity and associativity:
\begin{enumerate}
\item[(i)] $F(u,0)=F(0,u)=u$,
\item[(ii)] $F(u,v)=F(v,u)$,
\item[(iii)] $F(F(u,v),w)=F(u,F(v,w)$.
\end{enumerate}

The 
{\em Lazard ring}  
$\L$ is defined by the following construction \cite{Lazard}. 
Start with the polynomial ring 
$$\Z[\{A_{ij}, i,j\ge 1\}],$$ 
and form the power series 
\[
\tilde{F}(u,v):=u+v+\sum_{i,j\ge1}A_{ij}u^iv^j.
\]
Relation (i) is already satisfied.
Relations (ii) and (iii) 
give polynomial relations on the $A_{ij}$.  $\L$ is the quotient 
of $\Z[\{A_{ij}\}]$ by these relations. 
Letting $a_{ij}$ be the image of $A_{ij}$ in $\L$, 
the universal formal group law is 
$$F_\L(u,v)=u+v+\sum_{i,j\ge1}a_{ij}u^iv^j \in \L[[u,v]].$$
 We grade $\L$ by giving $a_{ij}$ degree $i+j-1$. 
If we give $u$ and $v$ degrees -1, then
has $F_\L(u,v)$ total degree $-1$.

To construct $\Omega_*$, we take the 
functor $\L_*\otimes_\Z\underline{\Omega}_*$ and impose the relations
\begin{multline*}
F_\L(\cn(L),\cn(M))(f:Y\to X,L_1,\ldots, L_r)\\
=\cn(L\otimes M)(f:Y\to X,L_1,\ldots, L_r)
\end{multline*}
for each pair of line bundles $L, M$ on $X$.

The construction of the pull-back in algebraic
cobordism    
for $\lci$ morphisms  is fairly technical, 
and is  the main task of \cite{ACII}.

The following universality statements are central
results  of \cite{ACI, ACII} (see \cite[Theorem 1.15]{ACII}).

\begin{thm}\label{thm:Universal} 
Algebraic cobordism is universal in both homology and cohomology:

\begin{enumerate}
\item[(i)]
 $X\mapsto \Omega_*(X)$
is the universal oriented Borel-Moore homology theory on $\Sch_k$. 
\item[(ii)]
 $X \mapsto \Omega^*(X)$ 
is the universal oriented cohomology theory on $\Sm_k$.
\end{enumerate}
\end{thm}

\noindent 
Let $A_*$ be an oriented Borel-Moore homology theory on $\Sch_k$.
Universality (i) yields a canonical natural transformation of
functors
$$\Omega_* \to A_*$$
which commutes with l.c.i pull-backs and external products.
Universality (ii) is parallel, see \cite{ACI}.
For the proof of Theorem \ref{thm:Universal}, the ground field
$k$ is required only to admit
resolution of singularities.

\section{Formal group laws}\label{sec:FGL} Let $A_*$ be
an oriented Borel-Moore homology theory on $\Sch_k$.
By \cite{ACI},
the Chern class of a tensor product is governed by a
formal group 
law $F_A(u,v)\in A_*(k)[[u,v]]$.
For each pair of line bundles $L, M$ on $Y\in\Sm_k$,
\begin{equation}\label{eqn:FGL0}
F_A(\cn(L),\cn(M))(1_Y)=\cn(L\otimes M)(1_Y).
\end{equation}
To make sense, the (commuting) operators $\cn(L)$ must
be nilpotent
on $1_Y\in A_*(Y)$. Nilpotency is proven in
 \cite {ACII}.

The existence of $F_A$, using the method employed by Quillen \cite{Quillen}, 
follows from an application 
of the projective bundle formula (PB) to a  product of 
projective spaces. We use the cohomological notation $A^*(\P^n\times \P^m)$.
By definition
\begin{align*}
A^*(\P^\infty\times\P^\infty)&=
\lim_{\substack{\leftarrow\\n,m}}A^*(\P^n\times\P^m)\\
&\cong \lim_{\substack{\leftarrow\\n,m}}A^*(k)[u,v]/(u^{n+1}, v^{m+1})\\
&=A^*(k)[[u,v]],
\end{align*}
where the isomorphism in the second line is
defined by sending 
$$au^iv^j \mapsto c_1(O_{\P^n\times\P^m}(1,0))^i
c_1(O_{\P^n\times\P^m}(0,1))^j
\cup p^*(a).$$
Here $a\in A_*(k)$, $p:\P^n\times\P^m\to\text{Spec}(k)$ is the 
structure morphism, and
\[
O_{\P^n\times\P^m}(i,j)=p_1^*O_{\P^n}(i)\otimes p_2^*O_{\P^m}(j).
\]

Clearly the 
elements 
$c_1(O_{\P^n\times\P^m}(1,1))\in A^1(\P^n\times\P^m)$ for varying $n, m$  
define an element $c_1(O(1,1))$ in the inverse limit. Therefore,
 there is a uniquely defined power series $F_A(u,v)\in A^*(k)[[u,v]]$ with
\[
c_1(O(1,1))=F_A(c_1(O(1,0),c_1(\sO(0,1)).
\]

If $Y\in\Sm_k$ is affine, then every pair of line bundles 
$$L, M \to X$$ is obtained by pull-back via a map 
$f:Y\to\P^n\times\P^m$ with 
$$L\cong f^*(O(1,0)), \ M\cong f^*(O(0,1)).$$
We conclude  
\[
c_1(L\otimes M)=F_A(c_1(L),c_1(M))
\]
by functoriality. Jouanolou's trick extends the equality
 to smooth quasi-pro\-jec\-tive $Y$.

For each oriented Borel-Moore homology 
theory $A_*$, there is a canonical graded ring homomorphism
\[
\phi_A:\L_*\to A_*(k)
\]
with $\phi_A(F_\L)=F_A$.

\begin{thm}[\hbox{\cite[Theorem 4]{ACI}}]\label{thm:Lazard} 
The homomorphism $\phi_\Omega:\L_*\to\Omega_*(k)$ is an isomorphism.
\end{thm}

Fix an embedding $\sigma:k\to \C$. Complex cobordism $MU^*(-)$ defines 
an oriented Borel-Moore cohomology theory $MU^{2*}_\sigma$ on $\Sm_k$ by
\[
X\mapsto MU^{2*}(X(\C)).
\]
By the universality of $\Omega^*$ as an oriented Borel-Moore 
cohomology theory on $\Sm_k$, we obtain a
 natural transformation $\vartheta^{MU,\sigma}:\Omega^*\to MU^{2*}_\sigma$. 
In particular, 
\[
\vartheta^{MU,\sigma}_\pt:\Omega^*(k)\to MU^{2*}(\pt).
\]
The formal group law for $MU^*$ is also 
the Lazard ring 
(after multiplying the degrees by 2, 
see \cite{Quillen}), so by Theorem~\ref{thm:Lazard}, 
the map $\vartheta^{MU,\sigma}_\pt$ is an isomorphism.

\label{qui}

\section{Oriented Borel-Moore functors of geometric type}
\label{gggg}
\subsection{Universality}
Algebraic cobordism
$\Omega_*$ is also a universal theory in the 
less structured setting of oriented Borel-Moore functors of 
geometric type. Since our goal will be to map $\Omega_*$ 
to the double point cobordism theory $\omega_*$, the less structure
required for $\omega_*$ the better.
We recall the definitions here.

\subsection{Oriented Borel-Moore functors with product}\label{obmp}
An {\em oriented Borel-Moore functor with product} on $\Sch_k$ consists
of
the following data:
\begin{enumerate}
\item[(D1)] An additive functor $H_*:\Sch_k' \to \Ab_*$.
\item[(D2)] For each smooth morphism
  $f:Y\to X$ in $\Sch_k$ of pure relative 
dimension $d$, a homomorphism of graded groups 
\[
f^*:H_*(X)\to H_{*+d}(Y).
\]
\item[(D3)] For each line bundle $L$ on $X$,
 a homomorphism of graded abelian groups
\[
\tilde{c}_1(L): H_*(X) \to H_{*-1}(X).
\]
\item[(D4)] 
For each pair $(X,Y)$ in $\Sch_k$, 
a bilinear graded pairing 
\begin{align*}
&\times :H_*(X) \times H_*(Y) \to H_*(X\times Y)\\
&(\alpha,\beta)\mapsto \alpha \times \beta
\end{align*} 
which is commutative, associative, and admits a
distinguished element $1\in H_0(\text{Spec}(k))$ as a unit.
\end{enumerate}
The pairing in (D4) is the {\em external product}.
The data (D1)-(D4) are required to satisfy eight conditions:
\begin{enumerate}
\item[(A1)]
Let  $f:Y\to X$ and $g:Z\to Y$ be smooth morphisms in $\Sch_k$
 of pure relative
dimension. Then, 
$$(f\circ g)^* = g^* \circ f^*.$$
Moreover, $\id_X^*=\id_{H_*(X)}$.
 
\item[(A2)] Let $f:X\to Z$ and $g:Y\to Z$ be morphisms 
in $\Sch_k$ where
$f$ is projective and $g$ is smooth of pure
relative dimension. In the cartesian square
\[
\xymatrix{
W\ar[r]^{g'}\ar[d]_{f'}&X\ar[d]^f\\
Y\ar[r]_g&Z\, ,}
\]
$f'$ is projective and $g$ is smooth or pure relative dimension. 
Then, 
\[
g^*f_*=f'_*g^{\prime*}.
\]
\item[(A3)] 
Let $f:Y\to X$ be projective. Then,  
\[
f_*\circ \cn(f^*L) =\cn(L)\circ f_*\, 
\]
for all line bundles $L$ on $X$.
\item[(A4)] Let
$f:Y\to X$  be smooth of pure relative dimension. Then, 
\[
\cn(f^*L)\circ f^* = f^*\circ \cn(L)\, .
\]
for all line bundles $L$ on $X$.
\item[(A5)] For all line bundles $L$ and $M$ on $X\in\Sch_k$,
\[
\cn(L)\circ \cn(M) = \cn(M) \circ \cn(L)\, .
\]
Moreover, if $L$ and $M$ are isomorphic, then
$\cn(L)=\cn(M)$. 
\item[(A6)] For projective morphisms $f$ and $g$,
\[
\times\circ (f_* \times g_*) = (f\times g)_*\circ\times\, .
\]
\item[(A7)] For smooth morphisms 
$f$ and $g$ or pure relative dimension,
\[
\times\circ(f^* \times g^*) = (f\times g)^*\circ\times\, .
\]
\item[(A8)] For $X,Y\in\Sch_k$,
\[
(\cn(L)(\alpha))\times \beta = \cn(p_1^*(L))\big(\alpha 
\times \beta\big),
\] 
for  $\alpha\in H_*(X)$, $\beta\in H_*(Y)$, and all line
bundles $L$ on $X$.

\end{enumerate}

 An oriented Borel-Moore homology theory $A_*$ on $\Sch_k$ determines an 
oriented Borel-Moore functor with product on $\Sch_k$ with 
the first Chern class operator is given by
\[
\cn(L)(\eta)=s^*s_*(\eta)
\]
for a a line bundle $L\to X$ with zero-section $s$.

 Let $H_*$ be an  oriented Borel-Moore functor 
with product on $\Sch_k$. The external 
products make $H_*(k)$ into a graded, commutative ring 
with unit $1\in H_0(k)$. For each $X$, 
the external product $$H_*(k)\otimes H_*(X)\to H_*(X)$$ makes 
$H_*(X)$  into a graded $H_*(k)$-module. 
The pull-back and push-forward maps are $H_*(k)$-module homomorphisms.

\subsection{Geometric type}\label{gtype}
Let $R_*$ be a graded commutative ring with unit. 
An {\em oriented Borel-Moore $R_*$-functor with product} 
on $\Sch_k$ is an  oriented Borel-Moore functor with product $H_*$ 
on $\Sch_k$ together with a graded ring homomorphism 
$$R_*\to H_*(k).$$

By the universal property of the Lazard ring $\L_*$, 
an oriented Borel-Moore $\L_*$-functor with product on $\Sch_k$ 
is the same as an oriented Borel-Moore functor with product $H_*$
 on $\Sch_k$ together with a formal group 
law $F(u,v)\in H_*(k)[[u,v]]$. In particular, 
an oriented Borel-Moore homology theory $A_*$ on $\Sch_k$ determines an 
oriented Borel-Moore $\L_*$-functor with product on $\Sch_k$.

An oriented Borel-Moore functor on $\Sch_k$ of {\em geometric type} 
is an oriented Borel-Moore $\L_*$-functor  $A_*$ with product on 
$\Sch_k$ satisfying the following
three additional axioms:
\begin{enumerate}
\item[(Dim)\,] 
 For $Y\in\Sm_k$ and line bundles  
$L_1,\dots,L_{r>\dim_k(Y)}$ on $Y$, 
\[
\cn(L_1)\circ\dots\circ\cn(L_r)(1_Y) = 0 \in A_*(Y)\, .
\]
\item[(Sect)\,]  
For $Y \in \Sm_k$ and
a section $s\in H^0(Y,L)$ of a line bundle $L$ transverse
to the zero section of $L$, 
\[
\cn(L)(1_Y) = i_*(1_Z),
\]
where $i:Z\to Y$ is the closed immersion of the zero subscheme  of $s$.
\item[(FGL)\,] 

For $Y\in\Sm_k$ and line bundles $L,M$ on  $Y$, 
\[
F_{A}(\cn(L),\cn(M))(1_Y) = \cn(L\oo M) (1_Y) \in A_*(Y)\, .
\]
\end{enumerate}
In axiom (FGL), $F_{A}\in A_*(k)[[u,v]]$ is the image of 
the power series $F_{\L}$ under
the homomorphism $\L_* \to A_*(k)$ giving the $\L_*$-structure. 

By 
\cite[Proposition 1.7]{ACII}, 
the oriented Borel-Moore functor with product on $\Sch_k$ 
determined by an oriented Borel-Moore homology theory on $\Sch_k$  
is an oriented Borel-Moore functor of geometric type.

\begin{thm}[\hbox{\cite[Theorem 2.17]{ACI}}] \label{thm:OBMGUniversal} 
The  oriented Borel-Moore functor of geometric type on $\Sch_k$ 
determined by $\Omega_*$ is universal.
\end{thm}

Let $A_*$ be an oriented Borel-Moore functor of
geometric type on $\Sch_k$.
Universality  yields a canonical natural transformation of
functors
$$\Omega_* \to A_*$$
which commutes with smooth pull-backs,
Chern class operators $\cn(L)$, and external products.
Again, only resolution of singularities for $k$ is
required for Theorem \ref{thm:OBMGUniversal}.

\section{The functor $\omega_*$}
\label{fff}
\subsection{Grading}
We start by defining a grading.
Let $X\in \Sch_k$. Let $$\sM_n(X)^+\subset \sM(X)^+$$
denote the subgroup generated by elements
$$[f:Y \to X]$$
where $Y\in \Sm_k$ is irreducible of dimension $n$.
Since the double point relations ${\mathcal R}(X)$
respect the dimension grading,
$$\omega_n = \sM_n(X)^+/{\mathcal R}_n(X)$$
defines a natural grading on $\omega_*(X)$.

\label{z}

\subsection{Push-forward, pull-back, and external products}
The 
assignment $X\mapsto \omega_*(X)$ carries the 
following elementary structures:\\

\noindent
{\em Projective push-forward}. Let $g:X\to X'$ be a 
projective morphism in $\Sch_k$. A map 
$$
g_*:\sM_*(X)^+\to \sM_*(X')^+
$$
is  defined by
$$g_*([f:Y\to X])= [g\circ f:Y\to X'].$$ 
By the definition of double point cobordism,  
$g_*$ descends to a functorial push-forward
$$g_*:\omega_*(X)\to \omega_*(X')$$
satisfying 
$$
(g_1\circ g_2)_*=g_{1*}\circ g_{2*}.$$

\vspace{10pt}
\noindent {\em Smooth pull-back}. Let $g:X'\to X$ be 
a smooth morphism in $\Sch_k$
of pure relative dimension $d$. 
A map 
$$g^*:\sM_*(X)^+\to \sM_{*+d}(X')^+$$ 
is defined by
$$g^*([f:Y\to X]) = [p_2:Y\times_X X'\to X'].$$
 Since the pull-back by $g\times\id_{\P^1}$ of a double point cobordism over $X$ is
a double point cobordism over $X'$, $g^*$ descends to a functorial pull-back 
 $$g^*:\omega_*(X)\to \omega_{*+d}(X')$$ 
satisfying
$$(g_1\circ g_2)^*=g_{2}^*\circ g_{1}^*.$$

\vspace{10pt}
\noindent
{\em External product}. A double 
point cobordism $\pi:Y\to X\times\P^1$ over $X$ 
gives rise to a double point cobordism 
$$Y\times Y'\to X\times X'\times\P^1$$ for each 
$[Y'\to X']\in \sM(X')$. Hence, the external product
\[
[f:Y\to X]\times [f':Y'\to X']=[f\times f':Y\times_k Y'\to X\times_k X']
\] 
on $\sM_*(-)^+$ descends  to an external product on $\omega_*$.

\vspace{10pt}
\noindent {\em Multiplicative unit}. The class 
$[\id:\text{Spec}(k) \to \text{Spec}(k)]\in \omega_0(k)$
is a unit for the external product on $\omega_*$.

\label{zz}

\subsection{Borel-Moore functors with product} \label{bmp}
A {\em Borel-Moore functor with product} on $\Sch_k$ consists
of the 
structures (D1), (D2), and (D4) of Section \ref{obmp}
satisfying axioms (A1),(A2), (A6),
and (A7). A Borel-Moore functor with product is
simply an oriented Borel-Moore functor with product {\em without}
Chern class operations.

\begin{lem}
Double point cobordism $\omega_*$ is a Borel-Moore functor with
product.
\end{lem}

\begin{proof}
The grading and
structures (D1), (D2), and (D4) have been constructed in
Sections \ref{z} and \ref{zz} .
Axioms (A1), (A2), (A6),
and (A7) follow easily from the definitions.
\end{proof}

\subsection{$\omega_*\to\Omega_*$} 
A natural transformation of Borel-Moore functors with
product is obtained once the double point 
relations are shown to be satisfied in $\Omega_*$.

Let $F(u,v)\in\Omega_*(k)[[u,v]]$ be the formal group law for $\Omega_*$. 
By definition,
 $$F(u,v)=u+v+\sum_{i,j\ge1}a_{i,j}u^iv^j$$
 with $a_{i,j}\in\Omega_{i+j-1}$. Let
$
F^{1,1}(u,v)=\sum_{i,j\ge1}a_{i,j}u^{i-1}v^{j-1}.
$
We have
\[
F(u,v)=u+v+uv\cdot F^{1,1}(u,v).
\]

Let $Y\in \Sm_k$.
Let $E_1$, $E_2$ be smooth divisors intersecting transversely in $Y$
with sum $E=E_1+E_2$. Let 
$$i_D: D=E_1 \cap E_2\to Y$$
be the inclusion of the intersection.
 Let
$O_D(E_1)$, $O_D(E_2)$ be the restrictions to $D$ of 
the line bundles $O_Y(E_1)$, $O_Y(E_2)$. 
Define an element $[E\to Y]\in \Omega_*(Y)$ by
\begin{multline*}
[E\to Y]=[E_1\to Y]+[E_2\to Y]\\
+i_{D*}\Big(F^{1,1}\big(\cn(O_D(E_1)),\cn(O_D(E_2))\big)(1_D)\Big)
\end{multline*}
The following result is proven in \cite{ACI} as
a consequence of the formal group law.

\begin{lem} \label{lem:DivFormula1}
Let $F\subset Y$ be a smooth divisor linearly equivalent
to $E$, then
\begin{equation*}
[F\to Y]=[E\to Y]\in \Omega_*(Y).
\end{equation*}
\end{lem}

\vspace{10pt}
If the additional condition
$$O_D(E_1)\cong O_D(E_2)^{-1}$$
is satisfied, a direct evaluation is possible.
 Let $\P_D\to D$ be the $\P^1$-bundle $\P(O_D \oplus O_D(E_1))$.

\begin{lem}\label{lem:DivFormula2}  We have
$$
F^{1,1}\big(\cn(O_D(E_1)),\cn(O_D(E_2))\big)(1_D)
= -[\P_D\to D] \in\Omega_*(D).$$
\end{lem}

\begin{proof} Both 
sides of the formula depend only upon 
the line bundles $O_D(E_1)$ and $O_D(E_2)$.
To prove the Lemma,  
we may replace $E$ with any $E'=E_1'+E_2'$ on any $Y'$, 
so long as $E_1'\cap E_2'=D$ and $O_{Y'}(E'_i)$ restricts to $O_D(E_i)$ on $D$.

The surjection $O_D\oplus O_D(E_1)\to O_D(E_1)$ defines a
section $s:D\to\P_D$ with normal bundle  $O_D(E_1)$. 
Let $Y'$ be the deformation to the normal cone 
of the closed immersion $s$. By definition,
$Y'$ is the blow-up of $\P_D\times\P^1$ along 
$s(D)\times0$.  The blow-up of $\P_D$ along $D$ is $\P_D$ and the
exceptional divisor $\P$ of $Y'\to P_D\times\P^1$ is also $\P_D$. 

The composition $Y'\to \P_D\times\P^1\to \P^1$ has
fiber $Y'_0$ over $0\in\P^1$ equal to $\P_D\cup \P$. 
The intersection $\P_D\cap \P$ is $s(D)$ and the line bundles
$O_{Y'}(\P)$, $O_{Y'}(\P_D)$ restrict to $ O_D(E_1)$, $O_D(E_2)$ on $s(D)$ 
respectively. Thus, we may use  $E'=Y'_0$, $E_1'=\P_D$, and $E_2'=\P$. 

By Lemma~\ref{lem:DivFormula1}, we have the relation 
$[Y'_\infty\to Y']=[Y'_0\to Y']$ in $\Omega_*(Y')$. By definition,
 $[Y'_0\to Y']$ is the sum
\begin{multline*}
[Y'_0\to Y]= [\P_D\to Y']+[\P\to Y']\\
+i_{D*}\Big(F^{1,1}\big(\cn(O_D(\P_D)),\cn(O_D(\P))\big)(1_D)\Big). 
\end{multline*}
Pushing forward  the relation $[Y'_\infty\to Y]
=[Y'_0\to Y']$ to $\Omega_*(D)$ by the 
composition $$Y'\to \P_D\times\P^1\xrightarrow{p_1} 
\P_D\to D$$ yields the relation
\[
[\P_D\to D]=[\P_D\to D]+[\P\to D]+F^{1,1}\big(
\cn(O_D(\P_D)),\cn(O_D(\P))\big)(1_D)
\] in $\Omega_*(D)$.
Since $\P\cong\P_D$ as a $D$-scheme, the proof is complete.
\end{proof}
\vspace{10pt}

\begin{cor}\label{dpl}
 Let $\pi:Y\to \P^1$ be a double point degeneration over
$0\in\P^1$. Let 
$$\pi^{-1}(0)=A\cup B.$$
Suppose the fiber $Y_\infty=\pi^{-1}(\infty)$ is smooth. Then
\[
[Y_\infty \to Y]=[A\to Y]+[B\to Y]-[\P(\pi)\to Y] \in \Omega_*(Y).
\]
\end{cor}

\vspace{10pt}
Sending $[f:Y\to X]\in \sM_n^+(X)$ to the class 
$[f:Y\to X]\in\Omega_n(X)$ defines a 
 natural transformation $\sM_*^+\to \Omega_*$ of 
Borel-Moore functors with product on $\Sch_k$. 

\begin{prop}\label{prop:Map} The map $\sM_*^+\to\Omega_*$ 
descends to a  natural transformation 
$$\vartheta:\omega_*\to \Omega_*$$ 
of Borel-Morel functors with product on $\Sch_k$. Moreover,
$\vartheta_X$ is surjective for each $X\in\Sch_k$.
\end{prop}

\begin{proof} Let $\pi:Y\to X\times\P^1$ be a double point degeneration
 over $X$. We obtain
a canonical
double point degeneration
 $$\pi'=(\id, p_2\circ \pi):Y\to Y\times\P^1.$$ 
Certainly
 $$\pi=(p_1\circ f,\id)\circ g.$$
 Since $\sM_*^+\to\Omega_*$ is compatible with projective push-forward, 
the first assertion reduces to Lemma \ref{dpl}.

The surjectivity follows from the fact that the canonical map 
\[
\sM_*(X)^+\to \Omega_*(X)
\]
 is surjective by \cite[Lemma 4.15]{ACI}.
\end{proof}

We will prove Theorem \ref{one} by showing $\vartheta$ is
an isomorphism.
The strategy of the proof is to show that $\omega_*$ admits first
Chern class operators and a
 formal group law and first Chern class operators, 
making $\omega_*$ into an oriented Borel-Moore functor of geometric type. 
We then use the universality of $\Omega_*$ given by 
Theorem~\ref{thm:OBMGUniversal} to determine
an inverse $\Omega_*\to \omega_*$ to $\vartheta$.

\section{Chern classes I}
\label{cc1}
Let $X\in \Sch_k$, and let 
 $L\to X$
be a line bundle generated by global sections.
 We will 
 define a {\em first  Chern class operator}
\[
\cn(L):\omega_*(X)\to \omega_{*-1}(X).
\]
A technical Lemma is required for the definition.

Let $[f:Y\to X]\in \sM(X)^+$ with
$Y$ irreducible of dimension $n$. 
For $s\in H^0(Y,f^*L)$, let $$i_s:H_s\to Y$$ be the inclusion of the 
zero subscheme of $s$. Let 
$$U\subset \P(H^0(Y,f^*L))=\{\ s \ | \ 
\text{$H_s$ is smooth and of codimension 1 in $Y$}\ \}.$$

\begin{lem}\label{lem:1stChern} 
We have
\begin{enumerate}
\item[(i)] $U$ is non-empty.
\item[(ii)] For $s_1, s_2\in U(k)$, 
$[H_{s_1}\to X]=[H_{s_2}\to X] \in \omega_{n-1}(X)$.
\end{enumerate}
\end{lem}

\begin{proof} Since $L$ is globally
generated, so is $f^*L$. Then (i) follows from Bertini's theorem
(using the characteristic 0 assumption for $k$).

Let
$\sH\subset  Y\times\P(H^0(Y,f^*L))$ be the universal
Cartier divisor.
 Let  $y\in Y$ be a closed point
with ideal sheaf $\mathfrak{m}_Y\subset \sO_Y$.
 Since $f^*L$ is globally generated, 
the fiber of $\sH\to Y$ over $y$ is the hyperplane 
$$\P(H^0(Y,f^*L\otimes \mathfrak{m}_y)\subset 
\P(H^0(Y,f^*L).$$ Hence, $\sH$ is smooth over $k$.

For (ii),
let $$i:\P^1\to  \P(H^0(Y,f^*L))$$ be a linearly embedded 
$\P^1$ with $i(0)=s_1$. By Bertini's theorem, 
the pull-back $$\sH_i=\sH\times_{\P(H^0(Y,f^*L))}\P^1$$ 
is smooth for general $i$. 
Clearly $\sH_i\to X\times\P^1$ gives a naive cobordism between 
$[H_{s_1}\to X]$ and 
$[H_{i(t)}\to X]$ for all $k$-valued points $t$ in a dense open subset of 
$\P^1$. Since $i$ is general, we have
 $$[H_{s_1}\to X]= [H_s\to X]\in \omega_{n-1}(X)$$
 for all  $k$-valued points $s$ in a dense open subset of $U$.
The same result for $s_2$ completes the proof.
\end{proof}

For $L$ globally generated, we can define the homomorphism
\[
\cn(L):\sM_*(X)^+\to\omega_{*-1}(X)
\]
by sending $[f:Y\to X]$ to $[H_s\to X]$ for 
 $H_s$ smooth and codimension 1 in $Y$.

\begin{lem}\label{lem:ChernDescent} The map $\cn(L)$ descends to
\[
\cn(L):\omega_*(X)\to\omega_{*-1}(X)
\]
\end{lem}

\begin{proof} Let $\pi:W\to X\times\P^1$ be a double point 
cobordism with degenerate fiber over $0\in \P^1$ and smooth fiber
over $\infty \in \P^1$. Hence,
$$W_0=A\cup B$$ with $A, B$ smooth divisors intersecting transversely 
in the double point locus $D=A\cap B$.
The double point relation is
\begin{equation}\label{eqn:relation}
[W_\infty \to X]= [A\to X]+[B\to X]-[\P(\pi)\to X].
\end{equation}

 Let $i_s:H_s\to W$ be the divisor of a general section $s$ of 
$(p_1\circ \pi)^*L$. As in the proof of lemma~\ref{lem:1stChern}, 
we may assume $H_s$, $H_s\cap W_\infty$, 
$H_s\cap S$, $H_s\cap A$ and $H_s\cap B$ 
are smooth divisors on $W$, $W_\infty$, $A$, $B$, and $D$
respectively. 
Then $$\pi\circ i_s:H_s\to X\times\P^1$$ is again a double point cobordism.
The associated double point relation 
\begin{equation*}
[H_s\cap W_\infty\to X]= 
[H_s\cap S\to X]+[H_s\cap T\to X]-[\P(\pi\circ i_s)\to X].
\end{equation*}
is obtained by applying
$\cn(L)$ term-wise to relation \eqref{eqn:relation}.
\end{proof}

Axioms (A3), (A4), (A5) and (A8) for
an oriented Borel-Moore functor with product are
easily checked for our definition of $\cn(L)$
{\em if all line bundles in question are globally generated}.
In particular, the operators $\cn(L)$ for globally
generated line bundles $L$ on $X$ are 
$\omega_*(k)$-linear and commute pairwise.

\begin{lem} \label{lem:dim} Let $X\in\Sch_k$, and let 
$$L_1,\ldots, L_{r>\dim_k X}\to X$$ be globally generated
 line bundles. Then, $$\prod_{i=1}^r\cn(L_i)=0$$ 
as an operator on $\omega_*(X)$.
\end{lem}

\begin{proof}  Let $[f:Y\to X]\in \sM(X)^+$. By
Bertini's theorem, $H_{f^*s}$ is smooth for a general choice 
of section $s\in H^0(X,L)$. Thus 
$$\cn(L)(f)=[f:H_{f^*s}\to X].$$
By induction, $\prod_i\cn(L_i)(f)$ is represented by the restriction of 
$f$ to $\cap_{i=1}^rH_{f^*s_i}$. But set-theoretically, 
$\cap_{i=1}^rH_{f^*s_i}= f^{-1}(\cap_{i=1}^rH_{s_i})$.
Since the sections $s_i$ are general, the intersection 
$\cap_{i=1}^rH_{s_i}$ is empty, whence the result.
\end{proof}

Let $F(u_1,\ldots, u_r)\in\omega_*(k)[[u_1,\ldots, u_r]]$ be a power series 
and let $L_1,\ldots, L_r$ be globally generated
on $X\in\Sch_k$. By Lemma~\ref{lem:dim}, the expression
$F(\cn(L_1),\ldots, \cn(L_r))$ is well defined
 as an operator on $\omega_*(X)$.

Lemma \ref{lem:dim} is condition (Dim) for an oriented
Borel-Moore functor of geometric type in case all the
line bundles in question are globally generated.

Chern classes for arbitrary line bundle will constructed in
Section \ref{cc2}. The axioms (FGL) and (Sect) will be verified
in Section \ref{cc2} and Section \ref{eee}.

\section{Extending the double point relation} 

\subsection{The blow-up relation} 
Before we construct the formal group law and the rest
of the Chern class operators for $\omega_*$, 
we describe two useful relations which are 
consequences of the basic double point cobordism relation.

The first is the blow-up relation.
Let $F\to X$ be a closed embedding 
in $\Sm_k$ with conormal bundle $\eta=\sI_F/\sI^2_F$ of rank $n$. Let 
$$\mu:X_F\to X$$ be the blow-up of $X$ along $F$. 
Let $\P_F$ be the $\P^{n-1}$-bundle $\P(\eta)\to F$. Let 
\begin{align*}
\P_1&=\P(\eta\oplus O_F)\to F\\
\P_2&=\P_{\P_F}(O_{\P_F}\oplus O(1)). 
\end{align*}
We consider $\P_1$ and $\P_2$ as $X$ schemes by the 
composition of the structure morphisms with the inclusion $F\to X$.

\begin{lem} \label{lem:blowup} We have
$$
[X_F\to X]=[\id:X\to X]-[\P_1\to X]+[\P_2\to X]\in\omega_*(X).$$
\end{lem}

\begin{proof} The Lemma
follows the
double point relation obtained
from the deformation to the normal cone of $F\to X$. 
Indeed, let 
$$\pi:Y\to X\times\P^1$$ be the 
blow-up along $F\times0$ with structure morphism 
$$\pi_2=p_2\circ\pi:Y\to \P^1.$$
 The fiber $\pi^{-1}(\infty)$ is just $X$,
 and 
$$\pi^{-1}(0) =X_F\cup \P_1,$$
 with $X_F$ and  $\P_1$ intersecting transversely along 
the exceptional divisor $\P_F$ of $\mu$. The normal 
bundle of $\P_F$ in $\P_1$ is $O(1)$. Thus
the associated double point relation is
\[
[\id:X\to X]= [X_F\to X]+[\P_1\to X]-[\P_2\to X]
\]
in $\omega_*(X)$.
\end{proof}

\subsection{The extended double point relation} 
Let $Y\in\Sm_k$. Let $A,B,C \subset Y$ 
be smooth divisors such that $A+B+C$ is 
a reduced strict normal crossing divisor. 
Let $$D=A\cap B, \ \ E=A\cap B\cap C.$$ 
As before, we let $O_D(A)$ denote the restriction of $O_Y(A)$ to $D$, 
and use a similar notation for the restrictions of bundles to $E$.
Let 
\begin{align*}
&\P_1=\P(O_D(A)\oplus O_D)\to D \\
&\P_E=\P(O_E(-B)\oplus O_E(-C))\to E\\
&\P_2=\P_{\P_E}(O\oplus O(1))\to\P_E\to E\\ 
&\P_3=\P(O_E(-B)\oplus O_E(-C)\oplus O_E)\to E. 
\end{align*}
We consider $\P_1$, $\P_2$ and 
$\P_3$ as $Y$-schemes by composing 
the structure morphisms with the inclusions $D\to Y$ and $E\to Y$.

\begin{lem}\label{lem:ExtendedDP} Suppose $C$ is linearly equivalent to $A+B$ 
on $Y$. Then,
\[
[C\to Y]=[A\to Y]+[B\to Y]-[\P_1\to Y]+[\P_2\to Y]-[\P_3\to Y]
\]
 in $\omega_*(Y)$
\end{lem}

\begin{proof} Let $Y_1\to Y$ be the blow-up of $Y$ along 
$(A\cup B)\cap C$.  
Since $(A\cup B) \cap C$ is a Cartier divisor on both $A\cup B$ and $C$,
 the proper transforms of both  $A\cup B$ and $C$
 define closed immersions $$A\cup B\to Y_1, \ \ C\to Y_1$$ 
lifting the inclusions $A\cup B\to Y$ and $C\to Y$.
 We denote the resulting closed subschemes of 
$Y_1$ by $A_1$, $B_1$ and $C_1$.

Let $f$ be a rational function on $Y$ with $\Div(f)=S+T-W$. We
obtain 
 a morphism $f:Y_1\to \P^1$ satisfying
$$f^{-1}(0)=A_1\cup B_1, \ \ f^{-1}(\infty)=C_1.$$ 
However, $Y_1$ is singular, unless $E=\0$. 
Indeed, if $A$, $B$ and $C$ are defined near a point $x$ of $E$ 
by local parameters $a$, $b$ and $c$, 
then locally analytically near $x\in A_1\cap B_1
\subset Y_1$, 
$$Y_1\cong E\times\Spec\left( k[a,b,c,z]/(ab-cz)\right).$$
Here, the exceptional divisor of $Y_1\to Y$ is defined by the ideal $(c)$,
   $A_1$ is defined by $(a,z)$ and $B_1$ is defined by $(b,z)$. 
The singular locus of $Y_1$ is isomorphic to $E$. 
We write $E_1$ for the singular locus of $Y_1$.

Let $\mu_2:Y_2\to Y_1$ be the blow-up of $Y_1$ along $A_1$. 
Since $A_1\subset Y_1$ is a Cartier divisor off of
 the singular locus $E_1$, the blow-up $\mu_2$ is an isomorphism
 over $Y_1\setminus E_1$. 
In our local description of $Y_1$, we see that $A_1\cap B_1$ 
is the Cartier divisor on $B_1$ defined by $(a)$, hence the 
proper transform of $B_1$ to $Y_2$ is isomorphic to $B$. Also, since
\[
b(a,z)=(ab,zb)=(zb, zc)=z(b,c),
\]
the strict transform of $A_1$ by 
$\mu_2$ is identified with the
blow-up $A_E$ of $A$ along $E$. In particular, 
since $E$ has codimension 2 in $A$ with normal bundle 
$O_E(B)\oplus O_E(C)$, we have the identification
\[
\mu_2^{-1}(E_1)=\P(O_E(-B)\oplus O_E(-C)).
\]
In addition, $Y_2$ is smooth. Indeed, 
the singular locus of $Y_2$ is contained in $$\mu_2^{-1}(E_1)
\subset\mu_2^{-1}(A_1)=A_E.$$ 
Since $A_E$ is a smooth Cartier divisor on $Y_2$, 
$Y_2$ is itself smooth, as claimed.

The morphism $\pi:Y_2\to \P^1$ defined by $\pi=f\circ \mu_2$ 
is a double point degeneration over
 $0\in \P^1$. with
$$\pi^{-1}(0)=A_E\cup B$$
 and  double point locus $A_E\cap B=A\cap B=D$.

Since $\pi^{-1}(\infty)=C$, we obtain the following
double point relation
\[
[C\to Y]=[A_E\to Y]+[B\to Y]- [\P(O_D(A)\oplus O_D)\to Y].
\]
in $\omega_*(Y)$.
Inserting the blow-up formula from Lemma~\ref{lem:blowup} completes the proof.
\end{proof}

 \section{Pull-backs in $\omega_*$} 
\subsection{Pull-backs}
The most difficult part of the construction of $\Omega_*$ is the 
extension of the pull-back maps from smooth morphisms to l.c.i. morphisms. 
We cannot hope to reproduce the full theory for $\omega_*$ directly.
Fortunately, only smooth pull-backs for $\omega_*$ are required
for the construction of an oriented Borel-Moore functor of geometric
type. However, our discussion of the formal group law for $\omega_*$
will require more than just smooth pull-backs.
The technique of {\em moving by translation} 
gives us sufficiently many pull-back maps for $\omega_*$.

\subsection{Moving by translation} 
\label{ccc}
We consider pull-back maps in the following setting. 
Let $G_1$ and $G_2$ be linear algebraic groups.
Let $Y\in \Sm_k$ admit a $G_1\times G_2$-action, and
let $B\in \Sm_k$ admit a transitive $G_2$-action.
Let
$$p:Y\to B$$
be a smooth morphism equivariant with respect to
$G_1\times G_2 \to G_2$.
Let 
$$s:B\to Y$$
be a section of $p$ satisfying three conditions:
\begin{enumerate}
\item[(i)] $s$ is
 equivariant with respect to
the inclusion $G_2\subset G_1\times G_2$,
\item[(ii)] $G_1\subset G_1\times G_2$ acts
trivially on $s(B)$,
\item[(iii)] $G_1 \times G_2$ acts transitively on $Y\setminus s(B)$.
\end{enumerate}
We will assume the above conditions hold throughout Section \ref{ccc}.

A special case in which all the
hypotheses are verified occurs
when $G_1=1$, $Y$ admits a
transitive $G_2$-action, and
$$p:Y \to Y, \ \ s: Y \to Y$$
are both the identity.

 \begin{lem}\label{lem:pullback}  
Let $i:Z \to Y$ be a morphism in $\Sm_k$ 
transverse to $s:B \to Y$. 
Let $f:W\to Y\times C$ be a projective morphism in $\Sm_k$.
 \begin{enumerate}
 \item[(1)]
 For all $g=(g_1,g_2)$ in a nonempty open set 
$$U(i,f)\subset G_1\times G_2,$$ 
the morphisms $(g\cdot i)\times\id_C$ and $f$ are transverse.
 \item [(2)] If $C=\Spec(k)$, then
for $g,g'\in U(i,f)$,
 \[
[Z\times_{g\cdot i} W\to Z] =  [Z\times_{g'\cdot i}W\to Z]\in \omega_*(Z).
 \]
 \end{enumerate}
 \end{lem}
 
 \begin{proof} Let $G=G_1\times G_2$. Consider the map 
$$\mu:G\times Z\to Y$$
defined by $\mu(g,z) = g \cdot i(z)$.
We first prove $\mu$ is smooth. In fact, we
will check $\mu$ is a submersion  
at each point $(g,z)$.

If $i(z)\in Y\setminus s(B)$, then $G\times z\to Y$ is 
smooth \footnote{Since $k$ has characteristic 0 and $G$ 
acts transitively on  $Y\setminus s(B)$, 
the orbit map is smooth.} 
and surjective by condition (iii), 
hence $\mu$ is a submersion at $(g,z)$ for all $g$. 

Suppose  $i(z)\in s(B)$. 
The map $G_2\times z\to s(B)$ is smooth and surjective by
condition (ii),
so the image of $T_{(g,z)}(G\times z)$ contains 
$$T_{i(z)}(s(B))\subset T_{i(z)}(Y).$$
Since $i$ is transverse to $s$,   $g\cdot i$ is transverse to $s$ for 
all $g$ and the composition
 \[
 T_zZ\xrightarrow{d(g\cdot i)}T_{g\cdot i(z)}(Y)\to 
T_{g\cdot i(z)}(Y)/T_{g\cdot i(z)}(s(B))
 \]
 is surjective. Thus 
 \[
 T_{(g,z)}(G\times Z)=T_{(g,z)}(G\times z)\oplus 
T_{(g,z)}(g\times Z)\xrightarrow{d\mu}
 T_{g\cdot i(z)}(Y)
 \]
 is surjective, and $\mu$ is a submersion at $(g,z)$.

The smoothness of $\mu$ clearly implies the smoothness of
$$\mu\times\id_C:G\times Z\times C\to Y\times C.$$
Hence $(G\times Z\times C)\times_\mu W$ is smooth over $k$,
 and the projection 
$$(G\times Z\times C)\times_\mu W\to
 G\times Z\times C$$
 is a well-defined element of $\sM(G\times Z\times C)$. 
Consider the projection
 \[
 \pi:(G\times Z\times C)\times_\mu W\to G.
 \]
 Since the characteristic is 0, 
the set of regular values of $\pi$ contains a nonempty
Zariski open dense subset 
$$U(i,f)\subset G.$$  Since
$G$ is an open subscheme of an affine space, 
the set of $k$-points of $U(i,f)$ is dense in $U(i,f)$. 
Any $k$-point $g=(g_1,g_2)$ in $U(i,f)$ satisfies claim (1) of
the Lemma.

For $g\in U(i,f)$, denote the element of $\sM(Z\times C)$ 
corresponding to 
$$(Z\times C) \times_{g\cdot i \times \id_C} W \to Z \times C$$
by $(g\cdot i)^*(f)$.
 
 For (2), let  $g, g'\in U(i,f)$ be two
$k$-points. We may consider $U(i,f)$ 
as an open subset of an affine space $\A^N$.
The pull-back $\pi^{-1}(\ell_{g,g'})$ of the 
line $\ell_{g,g'}$ through $g$ and $g'$ will be a closed subscheme of  
$(G\times Z)\times_\mu W$ which smooth and projective  over 
an open neighborhood $U\subset \ell_{g,g'}$ containing $g$ and $g'$. 
Then
\begin{equation} \label{ncd}
(U\times Z) \times_\mu W \to U
\end{equation}
provides a naive cobordism proving
\begin{equation}\label{edd}
[Z\times_{g\cdot i} W\to Z] =  [Z\times_{g'\cdot i}W\to Z]\in \omega_*(Z).
\end{equation}

Technically, the naive cobordism \eqref{ncd} has been
constructed only over
an open set $U \subset \P^1$. By taking a closure followed
by a resolution of singularities, the family \eqref{ncd} can be
extended appropriately over $\P^1$. The relation is \eqref{edd}
unaffected.
\end{proof}
 
Let $i:Z \to Y$ be a morphism in $\Sm_k$ of pure codimension $d$ transverse
to $s:B\to Y$.
We define
\begin{equation}\label{pbk}
i^*: \sM_*(Y)^+ \to \omega_{*-d}(Z)
\end{equation}
using (2) of Lemma \ref{lem:pullback} by
$$i^*[f:W \to Y] = [(g\cdot i)^*(f)]$$
for $g\in U(i,f)$.

 \begin{prop} \label{prop:pullback}
The pull-back  \eqref{pbk}
descends to a well-defined $\omega_*(k)$-linear pull-back 
 \[
 i^*:\omega_*(Y)\to \omega_{*-d}(Z).
 \]
 \end{prop}
 \begin{proof} 
The $\sM_*(k)^+$-linearity of the map
 $$i^*:\sM_*(Y)^+\to \omega_{*-d}(Z)$$ 
is evident from the construction.

Given a double point cobordism $f:W\to Y\times \P^1$ over
$0\in \P^1$, we will
show the pull-back of $f$ by  $(g\cdot i)\times\id_{\P^1}$ gives a 
double point cobordism for all $g$ in a dense open set of $U(i,f)$.

 Applying (1) of 
Lemma~\ref{lem:pullback} with $C=\P^1$ yields an open subscheme 
$$U_1\subset G_1\times G_2$$ for which $(g\cdot i)\times\id_{\P^1}$ pulls 
$W$ back to a smooth scheme $(g\cdot i)\times\id_{\P^1}(W)$, with a 
projective map to $Z\times \P^1$. Similarly, applying 
Lemma~\ref{lem:pullback} to the smooth fiber $W_\infty\to Y$, 
we find a subset  $U_2\subset U_1$ for which the 
fiber $(g\cdot i)\times\id_{\P^1}(W)_\infty$ is smooth. 
Finally, if $W_0=A\cup B$, applying 
Lemma~\ref{lem:pullback} to $A\to Y$, $B\to Y$ and $A\cap B\to Y$ 
yields an open subscheme $U_3\subset U_2$ for 
which $(g\cdot i)\times\id_{\P^1}(W)$ gives
the double point relation
\begin{multline*}
(g\cdot i)^*([W_\infty\to Y]) =\\
(g\cdot i)^*([A\to Y])+
(g\cdot i)^*([B\to Y])- (g\cdot i)^*([\P(f)\to Y]),
\end{multline*}
as desired.
\end{proof}

\begin{lem}\label{lem:ChernFunct} 
Let $L\to Y$ be a globally generated line bundle on $Y$. Then,
\[
i^*\circ\cn(L)=\cn(i^*L)\circ i^*.
\]
\end{lem}

\begin{proof} Since
 $i^*L$ is globally generated
on $Z$,  $\cn(i^*L)$ is well-defined. 
Let $[f:W\to Y]\in \sM(Y)$ and take $g\in G_1\times G_2$ so
 $g\cdot i:Z\to Y$ is transverse to $f$. 
For a general section $s$ of $f^*L$, the divisor of $s$, 
$$H_s\to W,$$
 is also transverse to $g\cdot i$. 
Hence,
$$i^*\circ\cn(L)([W\to Y]) = [Z \times_{g \cdot i} H_s \to Z].$$ 
Let
$H_{p_1^*(s)}$ be the divisor of $p_1^*(s)$ on $Z\times_{g\cdot i} W$
where $p_1$ is projection to the first factor.
Then,
$$\cn(i^*L)\circ i^*([W\to Y]) = [H_{p_1^*(s)} \to Z].$$ 
The isomorphism (as $Z$-schemes) 
$$Z\times_{g\cdot i}H_s\cong H_{p_1^*(s)}$$
yields the Lemma.
\end{proof}

\label{pbb}
\subsection{Examples}
There are two main applications of pull-backs constructed in Section \ref{pbb}.

First, let $Y=\prod_i \P^{N_i}$ be a product of projective spaces.
Let $$G_1=1, \  \ G_2=\prod_i \GL_{N_i+1}.$$
Let 
$p:Y \to Y$ and  $s: Y \to Y$
both be the identity.
For each morphism 
$$i:Z\to \prod_i\P^{N_i}$$ 
in $\Sm_k$ of codimension $d$, we have a well-defined 
$\omega_*(k)$-linear pull-back
\[
i^*:\omega_*(\prod_i\P^{N_i})\to \omega_{*-d}(Z).
\]

Second, let $Y$ be the total space of a line bundle $L$ on 
$B=\prod_i\P^{N_i}$ with projection $p$ and zero-section $s$,
$$p:L\to B, \ \ s:B\to L.$$
Here,
$G_1=\GL_1$ acts by scaling $L$, and
 $G_2=\prod_i\GL_{N_i+1}$ acts by symmetries on $B$.
For each morphism $$i:Z\to L$$ in $\Sm_k$
 which is transverse to the zero-section,
we have a $\omega_*(k)$-linear pull-back 
\[
i^*:\omega_*(L)\to \omega_{*-d}(Z).
\]

\label{exx}
\subsection{Independence}
\label{ind}
 The pull-backs constructed in
Section \ref{pbb} can be used to prove several independence
statements.

A {\em linear embedding} of $\P^{N-j} \to \P^N$ is an inclusion as
linear subspace. A {\em multilinear embedding}
$$\prod_{i=1}^m\P^{N_i-j_i}\to \prod_{i=1}^m\P^{N_i}$$
is a product of linear embeddings.
For fixed $j_i$, the multilinear embeddings are related
by naive cobordism. The classes
$$M_{j_1,\ldots,j_m} = \left[
 \prod_{i=1}^m\P^{N_i-j_i}\to \prod_{i=1}^m\P^{N_i}
\right] \in \omega_*(\prod_{i=1}^m\P^{N_i})$$
are therefore well-defined.

\begin{prop} \label{prop:ProductInd}The classes
$$\{ M_{j_1,\ldots,j_m}\ | \ 0\leq j_i \leq N_i \ \}
\subset \omega_*(\prod_{i=1}^m\P^{N_i})
$$  
are independent over $\omega_*(k)$.
\end{prop}

\begin{proof} 
Let $J=(j_1,\ldots, j_m)$ be a multi-index. There is a partial
ordering defined by
$$J\le J'$$
if $j_i\le j_i'$ for all $1\leq i \leq m$.
Let
$$\alpha=\sum_J a_J M_J \in \omega_*(\prod_{i=1}^m\P^{N_i})$$
where $a_J\in \omega_*(k)$.

If the $a_J$ are not all zero, let
$J_0=(j_1,\ldots, j_m)$ be a minimal multi-index
for which $a_J\neq0$. If we take a pull-back by a multi-linear embedding 
$$i:\prod_i\P^{j_i}\to \prod_i\P^{N_i},$$ then 
\[
i^*(\alpha)=a_{J_0}\cdot\left[
\prod_{i=1}^m\P^0\to \prod_{i=1}^m\P^{j_i}\right]
\in \omega_{*}(\prod_{i=1}^m\P^{j_i}).
\]
 Pushing-forward to $\omega_{*}(k)$ gives $a_{J_0}\neq 0$.
Hence $\alpha\neq 0$.
\end{proof}

Let $H_{n,m}\subset \P^n\times \P^m$ be the hypersurface defined 
by the vanishing of a general section of $O(1,1)$.
More generally, for $0\leq i\leq n$, 
let $$H^{(i)}_{n,m}\subset \P^n\times \P^m$$
 be the (smooth) subscheme defined by the 
vanishing of $i$ general sections of $O(1,1)$. 
Taking the linear embeddings $\P^{m-j}\to \P^n$, we may consider 
$$H^{(i)}_{n,m-j}\subset \P^n\times \P^m$$
for $0\le j\le m$. The proof of the following
result is identical to the proof of proposition~\ref{prop:ProductInd}.

\begin{lem}
The classes $[H^{(i)}_{n,m-j}\to \P^n\times\P^m]\in \omega_*(\P^n\times
\P^m)$ for
$0\le i\le n$, $0\le j\le m$ are independent over $\omega_*(k)$.
\end{lem}

If classes $H^{(i)}_{n,j}$ are taken for $i>n$,
we have a partial independence results.
 
  \begin{prop} \label{prop:ProductInd2}  If the identity
 \[
 \sum_{i=0}^{n+2m}
\sum_{j=0}^m\alpha_{i,j}\cdot [H^{(i)}_{n+m,m-j}\to \P^{n+m}\times\P^m]=0 
\in \omega_*(\P^{n+m}\times \P^m)
 \]
holds for $\alpha_{i,j} \in \omega_*(k)$,
then $\alpha_{i,j}=0$ for  $0\le i+j\le n+m$, $0\le j\le m$. 
\end{prop}

\begin{proof}  We argue by induction. Consider all pairs $(i,j)$ 
satisfying
$$0\le i+j\le n+m, \ \ 0\le j\le m$$
for which $\alpha_{i,j} \neq 0$.
 Of these, take the ones with minimal sum $i+j$, and of these, take the one with minimal $j$, denote the resulting pair by $(a,b)$. Note that $a\le a+b\le n+m$.

Take the pull-back of the identity
by a bi-linear embedding 
$$i:\P^a\times\P^b\to\P^{n+m}\times\P^m.$$
 Then, for each pair $(i,j)$ with $i+j>a+b$, 
$$i^*[H^{(i)}_{n+m,m-j}\to \P^{n+m}\times\P^m]=0,$$
since $H^{(i)}_{n+m,m-j}$ has codimension $i+j$. 
Similarly 
$$i^*[H^{(i)}_{n+m,m-j}\to \P^{n+m}\times\P^m]=0$$
 if $j>b$. Thus the identity in question pulls back to 
\[
\alpha_{a,b}\cdot [H^{(a)}_{a,0}\to \P^a\times\P^b] = 0
\]
Since $H^{(a)}_{a,0}=\Spec(k)$, pushing-forward to a point
yields $\alpha_{a,b}=0$.
\end{proof}

Let $Y_{N,M}$ be the total space of the bundle $O(1,-1)$ 
on $\P^N\times\P^M$, and let $Y_{i,j}\to Y_{N,M}$
be the closed immersion induced by the bi-linear 
embedding $$\P^i\times\P^j\to \P^n\times\P^m.$$

 \begin{prop} \label{prop:ProductInd3} 
If the identity
 \[
 \sum_{i=0}^N\sum_{j=0}^M \alpha_{i,j}\cdot [Y_{N-i,M-j}\to Y_{N,M}]=0 \in 
\omega_*(Y_{N,M})
 \]
holds for $\alpha_{i,j}\in\omega_*(k)$,
 then $\alpha_{i,j}=0$ for  $0\le i+j \le N$,  $0\le j\le M$. 
\end{prop}

\begin{proof} The proof is similar to that of 
Proposition~\ref{prop:ProductInd2}.
Consider all pairs $(i,j)$ 
satisfying
$$0\le i+j\le N, \ \ 0\le j\le m$$
for which $\alpha_{i,j} \neq 0$.
 Of these, take the ones with minimal sum $i+j$, 
and of these, take the one with minimal $j$, 
denote the resulting pair by $(a,b)$. Note that $a\le a+b\le N$.

Let $s_0,\ldots, s_N$ be  sections of
 $H^0(\P^a\times \P^b,\sO(1,1))$. 
Since $$N+1\ge a+b+1>\dim_k\P^a\times\P^b,$$
 we may choose the $s_i$ so as to have no common zeros. 
Hence  $s_0,\ldots, s_N$ define a morphism   
$$f:\P^a\times\P^b\to \P^N.$$
 Let $g:\P^b\to \P^M$ be a linear embedding. 
We obtain a morphism 
$$h=(f, g\circ p_2):\P^a\times\P^b\to \P^N\times\P^M$$ satisfying
$
h^*(O(1,-1))\cong O(1,0).$

A non-zero section $s\in H^0(\P^a\times\P^b,O(1,0))$ with smooth 
divisor defines a lifting $$(h, s):\P^a\times\P^b\to Y_{N,M}$$ of
$h$ which is 
transverse to the zero-section $$\P^N\times\P^M\to Y_{N,M}.$$

 We may therefore apply Proposition~\ref{prop:pullback} as explained in 
the second example of Section \ref{exx}
 to give a well-defined  $\omega_*(k)$-linear  pull-back map
\[
(h, s)^*:\omega_*(Y_{N,M})\to \omega_{*}(\P^a\times\P^b).
\]
We have 
\[
(h, s)^*([Y_{N-i,M-j}\to Y_{N,M}])=[H^{(i)}_{a,b-j}\to\P^a\times\P^b].
\]
Hence,
 $$(h, s)^*([Y_{N-i,M-j}\to Y_{N,M}])=0$$ 
 if $i+j>a+b$ or $j>b$ 
for dimensional reasons. Also,
\begin{eqnarray*}
(h, s)^*([Y_{N-a,M-b}\to Y_{N,M}])& = &[H^{(a)}_{a,0}\to\P^a\times\P^b]\\
& = &
[\Spec(k) \to \P^a\times\P^b].
\end{eqnarray*}
The pull-back of the identity 
stated in the Proposition by $(h,s)^*$ followed
by a push-forward to the point yields $\alpha_{a,b}=0$.
\end{proof}

\section{Admissible towers} 
\subsection{Overview}
We would like to construct a formal group 
law over $\omega_*(k)$ using  the method of Quillen described in 
Section~\ref{sec:FGL}. For
Quillen's construction,
the classes
\begin{equation}\label{clg}
\Big\{ [\P^i\times\P^j\to\P^n\times\P^m]\Big\}_
{0\le i\le n,\ 0\le j\le m} \subset \omega_*(\P^n \times \P^m)
\end{equation}
are required to constitute an $\omega_*(k)$-basis.
However, Lemma \ref{prop:ProductInd} only establishes independence.
We circumvent the problem by proving
a weak version of the generation of $\omega_*(\P^n\times\P^m)$ by the 
classes \eqref{clg}.

Let $Y$ be in $\Sm_k$.
An {\em admissible projective bundle} over $Y$
is a morphism of the form $$\P(\oplus_iL_i)\to Y$$
where the $L_i$ are line bundles on $Y$. An {\em admissible tower} 
over $Y$ 
is a morphism $\P\to Y$  which factorizes
\[
\P=\P_n\to \P_{n-1}\to\ldots\to \P_1\to \P_0=Y
\]
as a sequence of admissible projective bundles. The $i^{th}$
step 
$$\P_{i+1}\to \P_i$$
is an admissible projective bundle over $\P_i$.
We call $n$ the {\em length} of the admissible tower $\P\to Y$.
In particular, the identity $Y\to Y$ is an
admissible tower of length 0.

We prove the span of classes \eqref{clg} contains the
classes of all admissible towers over $\P^n\times \P^m$.

\subsection{Twisting}
Our main decomposition result for admissible towers
 $[\P\to Y]$ is based on twisting modifications in the
various steps of the tower. 

Let $Y\in \Sm_k$.
Let $E$ be a vector bundle on $Y$, let 
$L$ be a line bundle on $Y$,
and let $H$ a smooth divisor on $Y$. 
Let $E_H$, $L_H$ and $L(H)_H$ denote the restrictions to $H$. The projections
\begin{gather*}
E\oplus L\oplus L(H)\to E\oplus L,\\
E\oplus L\oplus L(H)\to E\oplus L(H)
\end{gather*}
give closed immersions 
\begin{gather*}
\P(E\oplus L)\to \P(E\oplus L\oplus L(H)),\\
\P(E\oplus L(H))\to \P(E\oplus L\oplus L(H)). 
\end{gather*}
The projective  bundle 
$$\P(E_H\oplus L_H\oplus L(H)_H)\to H$$
has a closed immersion over $H\to Y$, 
\[
\P(E_H\oplus L_H\oplus L(H)_H)\to  \P(E\oplus L\oplus L(H)). 
\]
The subvarieties $\P(E\oplus L)$,  $\P(E\oplus L(H))$, 
and $\P(E_H\oplus L_H\oplus L(H)_H)$ are smooth 
divisors in $\P(E\oplus L\oplus L(H))$. The union
$$\P(E\oplus L(H))+ \P(E\oplus L)+\P(E_H\oplus L_H\oplus L(H)_H)
\subset \P(E\oplus L\oplus L(H))
$$
has strict normal crossing singularities.

We also have the bundles 
\[
\P(E_H\oplus L_H)\to H,\  \ \P(E_H)\to H, 
\]
with closed immersions into $\P(E\oplus L\oplus L(H))$ over $H\to Y$.
The intersections
$$\P(E\oplus L)\cap \P(E_H\oplus L_H\oplus L(H)_H)=\P(E_H\oplus L_H),$$
$$\P(E\oplus L)\cap \P(E\oplus L(H)) \cap
\P(E_H\oplus L_H\oplus L(H)_H)
=\P(E_H)$$
are easily calculated.

\begin{lem}\label{lem:LinEquiv} 
The linear equivalence
\[
\P(E\oplus L(H))\sim \P(E\oplus L)+\P(E_H\oplus L_H\oplus L(H)_H)
\]
holds on $\P(E\oplus L\oplus L(H))$. 
\end{lem}

\begin{proof} 
Let $P$ denote $\P(E\oplus L \oplus L(H))$, and let 
$q: P          \to Y$ be the
structure morphism.
As $\P(E\oplus L)\subset P$ is given by the 
vanishing of the composition
\[
q^*(L(H))\to q^*(E\oplus L\oplus L(H))\to O_P(1),
\]
we find  $O_P(\P(E\oplus L))\cong q^*(L(H))^\vee(1)$. Similarly, 
\begin{gather*}
O_P(\P(E\oplus L(H)))\cong q^*(L)^\vee(1),\\
O_P(\P(E_H\oplus L_H\oplus L(H)_H))\cong q^*(O_Y(H)).
\end{gather*}
The linear equivalence of the Lemma is now easily
obtained.
\end{proof}

Let $H$ be a smooth divisor on $Y\in\Sm_k$. Let 
\[
\P=\P_n\to \P_{n-1}\to\ldots\to \P_1\to \P_0=Y
\]
be the factorization
 of an admissible tower $\P\to Y$ as a tower of admissible 
projective bundles. Fix an $i\le n-1$ and write the bundle $\P_{i+1}\to \P_i$ 
as $$\P(\oplus_{j=1}^rL_j)\to \P_i$$ for 
line bundles $L_j$ on $\P_i$. 

\begin{lem}\label{lem:TowerReduction}  
There exists an admissible tower $\P'\to Y$ which factors as
\[
\P'=\P'_n\to\P'_{n-1}\to\ldots\to\P'_{i+1}\to \P_i\to\ldots \P_1\to \P_0=Y
\]
with $\P'_{i+1}\to \P_i$ given by 
the bundle $$\P(\oplus_{j=1}^{r-1}L_j\oplus L_r(H))\to \P_i$$ and 
admissible towers $Q_0\to H$, $Q_1\to H$, $Q_2\to H$, $Q_3\to H$ satisfying
\[
[\P'\to Y]=[\P\to Y]+\sum_\ell(-1)^\ell i_{H*}([Q_i\to H]) \ \in
\omega_*(Y).
\]
\end{lem}

\begin{proof} If  $X\in \Sm_k$ is irreducible and
 $E\to X$  is a vector bundle, 
$$\Pic(\P(E))=\Pic(X)\oplus \Z\cdot [O(1)].$$ 
In particular, if $E\to F$ is a 
surjection of vector bundles on $X$, the restriction map 
$$\Pic(\P(E))\to \Pic(\P(F))$$ is surjective. 
Hence, if $\P_{\P(F)}\to \P(F)$ 
is an admissible projective bundle,
then there is an admissible projective bundle 
$\P_{\P(E)}\to \P(E)$ and an isomorphism of projective 
bundles over $\P(F)$ 
\[
\P_{\P(F)}\cong \P(F)\times_{\P(E)}\P_{\P(E)}.
\]
By induction on the length of an 
admissible tower, the same holds for each admissible tower $\P \to\P(F)$.

Let $E=\oplus_{i=1}^{r-1}L_i$, and let $L=L_r$.  Consider the
 admissible projective bundle 
$$\hat\P_{i+1}=\P(E\oplus L_r\oplus L_r(H))\to \P_i$$ 
and the closed immersions
\begin{align*}
&i_0:\P(E\oplus L)\to\hat\P_{i+1}\\
&i_1: \P(E\oplus L(H))\to \hat\P_{i+1}.
\end{align*}
By our remarks above, we may 
extend $i_0$ to a closed embedding of admissible towers over $Y$,
\[
\tilde{i}_0:\P\to \hat{\P},
\]
where $ \hat{\P}\to Y$ admits a factorization
\[
\hat{\P}=\hat{\P}_n\to \hat{\P}_{n-1}
\to\ldots\to\hat{\P}_{i+1}\to \P_i\to\ldots\to \P_1\to \P_0=Y
\]

Let $\tilde{i}_1:\P'\to \hat{\P}$ be the pull-back 
$\P(E\oplus L(H))\times_{\P_i}\hat{\P}$,
 and let $\hat{\P}_H\to H$ be the pull-back of 
$\hat{\P}\to Y$ via $H\to Y$. By Lemma \ref{lem:LinEquiv},
we have the linear equivalence 
\[
\P'\sim  \P+\hat{\P}_H
\]
on the admissible tower $\hat{\P}$.

Since $\P(E\oplus L)+\P(E_H\oplus L_H\oplus L(H)_H)+\P(E\oplus L(H))$ 
is a reduced  strict normal crossing divisor on 
$\P(E\oplus L\oplus L(H))$, the sum  
$\P+\hat{\P}_H+\P'$ is a reduced strict normal crossing divisor 
on $\hat{\P}$. Since 
\begin{align*}
&\P(E\oplus L) \cap \P(E_H\oplus L_H\oplus L(H)_H)=\P(E_H\oplus L_H),\\
&\P(E\oplus L)\cap 
\P(E_H\oplus L_H\oplus L(H)_H)\cap \P(E\oplus L(H))=
\P(E_H)\end{align*}
are both admissible projective bundles over $\P_i\times_YH$, 
$$D=\P \cap \hat{\P}_H,\ \ F=\P \cap
\hat{\P}_H \cap \P'$$
 are both admissible towers over $H$. 
Let
\begin{align*}
&Q_0=\hat{\P}_H\\
&Q_1=\P_D(O_D(\P)\oplus O_D)\\
&Q_2=\P_{\P_F(O_F(-H)\oplus O_F(-\P'))}(O\oplus O(1))\\
&Q_3=\P_F(O_F(-H)\oplus O_F(-\P')\oplus O_F).
\end{align*}
Each $Q_i\to H$ is an admissible tower. 
Lemma~\ref{lem:ExtendedDP} completes the proof.
\end{proof}

\subsection{Generation}
Let $\omega_*(k)'\subset \omega_*(k)$ be the subgroup 
generated by  classes of admissible towers 
over $\Spec(k)$. Clearly,  $\omega_*(k)'$ is a subring.
 
Let $H_1,\ldots, H_s$ be divisors on $Y\in \Sm_k$ 
for which the associated invertible sheaves $\sO_Y(H_i)$ 
are generated by global sections. 
Let $$I=(i_1,\ldots, i_s)$$ be a multi-index with $i_r$ non-negative
for all $r$.
Let $$[H^I\to Y]\in\omega_*(Y)$$ denote the class
 of the closed immersion $H^I\to Y$, 
where $H^I$ is the closed subscheme of codimension
$\sum_r i_r$ defined by the
simultaneous vanishing of 
$i_1$ sections of $\sO_Y(H_1)$,
$i_2$ sections of $\sO_Y(H_2)$, \dots, and
$i_s$ sections of $\sO_Y(H_2)$.
By definition,
$$[H^{(0,\ldots,0)}\to Y] = [Y\to Y].$$
For a general choice of sections, $H^I$ is smooth. 
By naive cobordisms, $[H^I\to Y]$ is independent of the choice of sections.

The subvarieties $H^I$ may not be irreducible.
Let $H^I_1,\ldots H^I_{n_I}$ be the irreducible components of $H^I$.

\begin{lem} \label{lem:TowerDecomp} 
If
the restrictions
of the invertible sheaves $\sO_Y(H_i)$ generate
$\Pic(H^I_j)$ for every $H^I_j$,
then the classes of admissible towers over $Y$
lie in the $\omega_*(k)'$-span of $[H^I_j\to Y]$ in
$\omega(Y)$.
\end{lem}

\begin{proof} Given an admissible tower $\P\to Y$, we must
find an identity
\[
[\P\to Y]=\sum_{I,j} a_{I,j}\cdot [H^I_j\to Y]\ \in \omega(Y)
\]
with $a_{I,j}\in \omega_*(k)'$.

 We may assume $Y$ is irreducible and the divisors $H_i$ are smooth. 
If $Y$ has dimension 0, then every line bundle on $Y$ is trivial. 
By induction on the length of the tower,
every admissible tower $\P\to Y$ is the pull-back of an admissible 
tower $\P'\to \Spec(k)$ by the structure morphism 
$Y\to\Spec(k)$. The result is proven
in case $\dim_k Y =0$.

We proceed by induction on $\dim_k Y$. Let $\omega_*(Y)'$ be the  
subgroup generated by the push-forward to $Y$ of classes of the 
form $[\P'\to H^I_j]$, 
where $\P'\to H^I_j$ is an admissible tower and 
$I\neq(0,\ldots, 0)$. 
Since such $H^I_j$ satisfy the hypotheses of the
Lemma and have dimension strictly less than $Y$,
the push-forwards to $Y$ of the 
classes $[\P'\to H^I_j]$ lie in the $\omega_*(k)'$-span of
the classes $[H^I_j \to Y]$.

Let $\P \to Y$ be an admissible tower of length $n$
which factors
as
\[
\P\to Q\to Y
\]
where $\P\to Q$ is an admissible tower of length $n-i$ and
 $Q\to Y$ is an admissible tower of length $i<n$
isomorphic to a pull-back
$$Q\cong Q_0\times_k Y\to Y$$
of an admissible tower 
$Q_0\to\Spec(k)$ of length $i$.
By twisting, 
we will prove the condition
\begin{equation}\label{clm}
[\P\to Y]-[\P'\to Y] \in \omega_*(Y)'
\end{equation}
is satisfied for an
admissible tower $\P'\to Y$ of length $n$ which
 admits a factorization $\P'\to Q'\to Y$ as above where
 $Q'\to Y$ is an admissible tower of length $i+1$ of the 
form $$Q'\cong Q_0'\times_k Y\to Y$$ 
for an admissible tower $Q'_0\to \Spec(k)$ of length $i+1$.

The construction of $\P'\to Y$ satisfying \eqref{clm}
follows directly from Lemma~\ref{lem:TowerReduction}. 
Indeed, suppose 
$$\P_{i+1}\to \P_i=Q$$ 
is of the form $\P_Q(\oplus_i L_i)\to Q$. Since $Q=Q_0\times_kY$, we have
\[
\Pic(Q)=\Pic(Q_0)\oplus \Pic(Y).
\]
We can write each $L_i$ as
\[
L_i\cong p_1^*L_i^0\otimes p_2^* M_i
\]
for suitable line bundles $L_i^0$ on $Q_0$, and $M_i$ on $Y$.  
By Lemma~\ref{lem:TowerReduction}, 
the class $[\P\to Y]$ is equivalent modulo $\omega_*(Y)'$ to a 
class $[\tilde{\P}\to Y]$, where $\tilde{\P} \to Y$ is an admissible tower of 
length $n$ which factors as
\[
\tilde{\P}\to \tilde{\P}_{i+1}\to Q\to Y
\]
and where 
$\tilde{\P}_{i+1}=\P(\oplus_{i\neq j}L_i\oplus L_j(H_\ell))$ 
for any choice of $j$ and $\ell$ we like. 
Since the $H_\ell$ generate $\Pic(Y)$, after several such
applications of  Lemma~\ref{lem:TowerReduction}, we may
replace $\P$ with an admissible tower 
$$\P' \to \P'_{i+1}\to Q\to Y,$$ where
\[
\P'_{i+1}\cong \P(\oplus_i\ p_1^*L_i^0\otimes p_2^*L)\cong 
\P(\oplus_i\ p_1^*L_i^0)
\]
for a line bundle $L$ on $Y$. 
Thus ${\P'}_{i+1}\to Q\to Y$ is the pullback to $Y$ 
of an admissible tower $Q'_0\to Q_0\to \Spec(k)$,
and we obtain condition \eqref{clm}.

Repeated application of \eqref{clm} yields the relation
$$[\P\to Y] - [Q \to Y] \ \in \omega_*(Y)'$$
where 
$$Q\cong Y \times _k Q_0 \to Y$$
for an admissible tower $Q_0\to \Spec(k)$ of length $n$.
\end{proof}

\begin{cor} \label{prop:decomp1}
Let $\P\to \prod_{i=1}^m\P^{N_i}$ be an admissible tower. 
Then,
\[
[\P\to  \prod_{i=1}^m\P^{N_i}]= \sum_{J=(j_1,\ldots, j_m)}
a_J\cdot M_j \ \in \omega_*( \prod_{i=1}^m \P^{N_i})\]
for unique elements $a_J\in\omega_*(k)'$.
\end{cor}

\begin{proof}
For existence, we apply Lemma~\ref{lem:TowerDecomp} 
with $Y=\prod_{i=1}^m\P^{N_i}$ and  the divisors
$H_i$ defined by the pull-backs of 
hyperplanes in $\P^{N_i}$ via the projections $Y\to \P^{N_i}$. 
Uniqueness  follows from Proposition~\ref{prop:ProductInd}.
\end{proof}

\begin{cor}\label{prop:decomp11} 
Let $\P\to H_{n,m}$ be an admissible tower. 
Then, 
\[
[\P\to H_{n,m}]=\sum_{i,j}a_{i,j}\cdot[H_{n-i,m-j}\to H_{n,m}]
\]
for elements $a_{i,j}\in\omega_*(k)'$.
\end{cor}

Here, $H_{n-i,m-j}\to H_{n,m}$ is induced by the bi-linear embedding 
$$\P^{n-i}\times\P^{m-j}\to \P^n\times\P^m.$$
The sum in Corollary \ref{prop:decomp11} is over 
$$0\leq i \leq n, \ \ 0\leq j \leq m, \ \ i+j<n+m$$
for dimension reasons.

\begin{proof} 
We apply Lemma~\ref{lem:TowerDecomp} with 
$Y=H_{n,m}$ and divisors
$H_1=H_{n-1,m}$, $H_2=H_{n,m-1}$.
If $n\ge m$, the projection 
$$p_2:H_{n,m}\to\P^m$$ 
expresses  $H_{n,m}$ as a $\P^{n-1}$-bundle over $\P^m$.
Hence, $H_1$ and $H_2$ generate $\Pic(H_{n,m})$. Since 
\[
H_1^{(i)}\cdot H_2^{(j)}=H_{n-i,m-j},
\]
the hypotheses of Lemma~\ref{lem:TowerDecomp} 
are satisfied and yield the desired result.
\end{proof}

\begin{prop} \label{prop:decomp2} Let  $\P\to H_{n,m}$ be an admissible tower.
 Then,
\[
i_{H_{n,m}*}([\P\to H_{n,m}])=
\sum_{(i,j)\neq (0,0)}a_{i,j}\cdot [\P^{n-i}\times\P^{m-j}\to \P^n\times\P^m]
\]
for unique elements $a_{i,j}\in\omega_*(k)'$. 
\end{prop}

\begin{proof} 
If  $m=0$, then $H_{n,m}$ is a hyperplane in $\P^n$, and the result follows 
from Corollary~\ref{prop:decomp1}.
The same argument is valid for  $n=0$. 

We proceed by induction on $(n,m)$. Only
existence is required since
uniqueness follows from Proposition~\ref{prop:ProductInd}.
By Corollary~\ref{prop:decomp11}, we need only construct 
relation of the form
\[
i_{H_{n,m}*}(a\cdot [H_{n,m} \to H_{n,m}])=\sum_{(i,j)\neq (0,0)}
a_{i,j}\cdot [\P^{n-i}\times\P^{m-j}\to \P^n\times\P^m].
\]
for $a_{i,j}\in\omega_*(k)'$ for every $a\in \omega_*(k)$.
Since $$i_{H_{n,m}*}(a\cdot [H_{n,m}\to H_{n,m}])
=a\cdot i_{H_{n,m}*}([H_{n,m}\to H_{n,m}]),$$
the case $a=1$ suffices.

We have the linear equivalence on $\P^n\times\P^m$, 
\[
H_{n,m}\sim \P^{n-1}\times\P^m+\P^n\times\P^{m-1}.
\]
By the extended double point relation of Lemma~\ref{lem:ExtendedDP}, 
there are admissible towers $\P_1\to \P^{n-1}\times\P^{m-1}$, 
$\P_2\to H_{n-1,m-1}$ and $\P_3\to H_{n-1,m-1}$ for which
\begin{eqnarray*}
[H_{n,m}\to \P^n\times\P^m]&=& \ \ [\P^{n-1}\times\P^m\to \P^n\times\P^m]\\
& & + [\P^{n}\times\P^{m-1}\to \P^n\times\P^m]\\
& & -[\P_1\to\P^n\times\P^m]\\
& & +[\P_2\to \P^n\times\P^m]\\
& & -[\P_3\to \P^n\times\P^m].
\end{eqnarray*}
By induction, the classes $[\P_2\to \P^{n-1}\times\P^{m-1}]$ 
and $[\P_3\to \P^{n-1}\times\P^{m-1}]$ are expressible as
\[
[\P_\ell\to \P^{n-1}\times\P^{m-1}]=
\sum_{i,j}a^\ell_{i,j}\cdot
[\P^{n-i-1}\times\P^{m-j-1}\to \P^{n-1}\times\P^{m-1}],
\]
for $a^\ell_{i,j} \in \omega_*(k)'$ for
$\ell=2,3$. By Corollary~\ref{prop:decomp1}, a similar expression 
is obtained in case $\ell=1$.
\end{proof}

\section{The formal group law over $\omega_*(k)$}  
\label{rrr}
We use the classical method of Quillen to construct a 
formal group law over $\omega_*(k)$. 
Proposition~\ref{prop:decomp2} replaces the projective bundle formula.

By Proposition~\ref{prop:decomp2}, there are unique 
elements $a_{i,j}^{n,m}\in\omega_{i+j-1}(k)$
for which the identity
\begin{equation}\label{eqn:FGL1}
[H_{n,m}\to\P^n\times\P^m]=\sum_{(i,j)\neq (0,0)}
a_{i,j}^{n,m}\cdot [\P^{n-i}\times\P^{m-j}
\to\P^n\times\P^m] 
\end{equation}
holds in $\omega_*(\P^n\times \P^m)$.
For convenience, we set $a_{0,0}^{n,m}=0$.

\begin{lem} \label{stz}
If $N\ge n$, $M\ge m$, then $$a_{i,j}^{N,M}=a_{i,j}^{n,m}$$ 
for $0\le i\le n$, $0\le j\le m$.
\end{lem} 

\begin{proof} Pull-back 
relation \eqref{eqn:FGL1} for $N,M$ by a bi-linear 
embedding $$i:\P^n\times\P^m\to \P^N\times\P^M,$$
see Section \ref{exx} for the pull-back construction. 
We find
\begin{align*}
&i^*([H_{N,M}\to\P^N\times\P^M])=[H_{n,m}\to\P^n\times\P^m]\\
&i^*([\P^{N-i}\times\P^{M-j}\to\P^N\times\P^M])=
[\P^{n-i}\times\P^{m-j}\to\P^n\times\P^m]
\end{align*}
for $0\le i\le n$ and $0\le j\le m$. Since
 $i^*$ is $\omega_*(k)$-linear, 
the result follows from the uniqueness of the $a_{i,j}^{n,m}$.
\end{proof}

By Lemma \ref{stz}, we may define $a_{i,j}\in \omega_*(k)$ by
\[
a_{i,j}=\lim_{N\to\infty, M\to \infty}a_{i,j}^{N,M}.
\]
Following the convention 
$$[\P^{n-i}\times\P^{m-j}\to\P^n\times\P^m]=0$$
 if $i>n$ or if $j>m$, we write $a_{i,j}$ for $a_{i,j}^{n,m}$ 
in relation \eqref{eqn:FGL1}.

Taking $n=0$ and noting $H_{0,m}=\P^{m-1}$ 
 linearly embeds in $\P^m$, we find
\[
a_{0,1}=1, \ \ a_{0,j>1}=0.
\]
As the exchange of 
factors $\P^n\times\P^m\to \P^m\times\P^n$ sends 
$H_{n,m}$ to $H_{m,n}$, we obtain the symmetry
\[
a_{i,j}=a_{j,i}.
\]
Let $F_\omega(u,v)\in\omega_*(k)[[u,v]]$ be the power series
\[
F_\omega(u,v)=u+v+\sum_{i,j\ge1}a_{i,j}u^iv^j.
\]

\begin{prop}\label{prop:FGL2} 
Let $L_1$ and $L_2$ be globally generated
 line bundles on  $X\in\Sch_k$. Then, $L_1\otimes L_2$ is globally
generated and
\[
\cn(L_1\otimes L_2)=F_\omega(\cn(L_1),\cn(L_2)).
\]
\end{prop}

\begin{proof} The Lemma follows from the equation
\begin{equation}\label{ccm}
\cn(L_1\otimes L_2)(1_Y)=F_\omega(\cn(L_1),\cn(L_2))(1_Y).
\end{equation}
for all $L_1, L_2$ on all $Y\in\Sm_k$.
Indeed, if $[f:Y\to X]\in\sM(X)^+$, 
then 
$$f_*(1_Y)=[f:Y\to X] \in\omega_*(X).$$ 
By (A3), we have
\[
\cn(L)([f:Y\to X])=\cn(f_*(1_Y))=f_*(\cn(f^*L)(1_Y))
\]
for all globally generated $L$ on $X$, which verifies the claim.

Since $L_1$ and $L_2$ are globally generated, 
we have morphisms $$f_i:Y\to \P^{n_i}$$ with 
$L_i\cong f_i^*(O(1))$
for $i=1,2$. Thus,
 $$L_1\otimes L_2\cong (f_1\times f_2)^*(O(1,1)).$$

By the functoriality of 
Lemma~\ref{lem:ChernFunct},
we need only  prove \eqref{ccm} in  case 
$$Y=\P^n\times\P^m,\  L_1=O(1,0),\ L_2=O(0,1),\  L_1\otimes L_2=O(1,1).$$
Since
\begin{align*}
&\cn(O(1,1))(1_{\P^n\times\P^m})=[H_{n,m}\to \P^n\times\P^m]\\
&\cn(O(1,0))^i\circ \cn(O(0,1))^j(1_{\P^n\times\P^m})=
[\P^{n-i}\times\P^{m-j}\to \P^n\times\P^m],
\end{align*}
the defining  relation \eqref{eqn:FGL1} for the $a_{i,j}$ becomes
\[
\cn(O(1,1))(1_{\P^n\times\P^m})=
F_\omega(\cn(O(1,0),\cn(O(0,1))(1_{\P^n\times\P^m}),
\]
as desired.
\end{proof}

\begin{prop} \label{formgp}
$F_\omega(u,v)$ defines a formal group law over $\omega_*(k)$.
\end{prop}

\begin{proof} Of the axioms for formal group laws, the
first two have already been established:
\begin{enumerate}
\item[(i)] $F(u,0)=F(0,u)=u$,
\item[(ii)] $F(u,v)=F(v,u)$.
\end{enumerate}
The last axiom
\begin{enumerate}
\item[(iii)] $F(F(u,v),w)=F(u,F(v,w))$.
\end{enumerate}
will now be proven.

Let $G_1(u,v,w)=F(F(u,v),w)$ and $G_2(u,v,w)=F(u,F(v,w))$. 
For $\ell=1,2$, write
\[
G_\ell(u,v,w)=\sum_{i,j,k}a^\ell_{i,j,k}u^iv^jw^k.
\]
For globally generated
 line bundles $L_1, L_2, L_3$ on $X\in \Sch_k$,
\[
G_1(\cn(L_1),\cn(L_2),\cn(L_3))=F(\cn(L_1\otimes L_2),\cn(L_3))=\cn(L_1\otimes L_2\otimes L_3)
\]
by Proposition~\ref{prop:FGL2}.
A similar equation holds for $G_2$. Thus
\[
G_1(\cn(L_1),\cn(L_2),\cn(L_3))=G_2(\cn(L_1),\cn(L_2),\cn(L_3))
\]
as operators on $\omega_*(X)$. 

Specializing to  $X=\P^n\times\P^m\times \P^r$, we find
\begin{multline*}
G_\ell(\cn(O(1,0,0),\cn(O(0,1,0)),\cn(O(0,0,1))(1_X)\\
=\sum_{i=0}^n\sum_{j=0}^m\sum_{k=0}^ra^\ell_{i,j,k}\cdot
[\P^{n-i}\times\P^{m-j}\times\P^{r-k}\to \P^n\times\P^m\times\P^r]
\end{multline*}
for $\ell=1,2$. By Proposition \ref{prop:ProductInd}, 
\[
a_{i,j,k}^1=a_{i,j,k}^2
\]
for $0\le i\le n$, $0\le j\le m$, $0\le k\le r$. 
As $n$, $m$ and $r$ were arbitrary, the proof is complete.
\end{proof}

\section{Chern classes II}
\label{cc2}
\subsection{Definition}
Because $F_\omega(u,v)$ is a formal group law, there exists
an  inverse power series $\chi_\omega(u)\in\omega_*(k)[[u]]$
 characterized by the identity
\[
F_\omega(u,\chi_\omega(u))=0.
\]
We let $F^-_\omega(u,v)$ be the difference in our group law, 
\[
F^-_\omega(u,v)=F_\omega(u,\chi_\omega(v)).
\]

Using $F^-_\omega(u,v)$, we can extend the definition of $\cn(L)$ given
in Section \ref{cc1} for
globally generated $L$ to arbitrary line bundles. 

\begin{lem}\label{lem:FirstChernExt} 
Let $L, M, N$ be line bundles on $Y\in \Sm_k$ where
 $$L,\ M,\ L\otimes N,\ M\otimes N$$ are globally generated. Then,
\[
F^-_\omega(\cn(L),\cn(M))=F^-_\omega(\cn(L\otimes N), \cn(M\otimes N))
\]
as operators on $\omega_*(Y)$.
\end{lem}

\begin{proof} We first assume $N$ is globally generated. Then 
\begin{eqnarray*}
\cn(L\otimes N)&=&F_\omega(\cn(L),\cn(N))\\
\cn(M\otimes N)&=&F_\omega(\cn(M),\cn(N))
\end{eqnarray*} 
by Proposition \ref{prop:FGL2}.
The result then follows from the power series identity
\[
F^-_\omega(F_\omega(u,w),F_\omega(v,w))=F^-_\omega(u,v).
\]

In general, since $Y$ is quasi-projective, 
there is a very ample line bundle $N'$ such that $N''=N'\otimes N^{-1}$ 
is very ample. Then
\begin{eqnarray*}
F^-_\omega(\cn(L),\cn(M)) &=&F^-_\omega(\cn(L\otimes N'),\cn(M\otimes N'))\\
&=&
F^-_\omega(\cn(L\otimes N\otimes N''),
\cn(M\otimes N\otimes N''))\\
&=& F^-_\omega(\cn(L\otimes N),\cn(M\otimes N)),
\end{eqnarray*}
completing the proof. \end{proof}

Let $L$ be an arbitrary 
line bundle on $X\in\Sch_k$. Define the operator
\[
\cn(L):\sM_*(X)^+\to\omega_{*-1}(X)
\]
by the  following construction.
Let $Y\in \Sm_k$ be irreducible.
Let
\begin{equation}\label{ggg} 
[f:Y\to X]\in \sM(X)^+.
\end{equation}
Let $M$ be a very ample line bundle on $Y$ for which
$f^*(L) \otimes M$ is also very ample. 
Then,
\[
\cn(L)([f:Y\to X])=f_*\Big(F^-_\omega\big(
\cn(f^*(L)\otimes M),\cn(M)\big)(1_Y)\Big).
\]
By Lemma~\ref{lem:FirstChernExt},
 $\cn(L)([f])$ is  independent of the choice of $M$.
 Since $\sM_*(X)^+$ is the free abelian group with
generators \eqref{ggg},  $\cn(L)$ is defined on $\sM_*(X)^+$.

Let $X\in \Sch_k$, and
let $\pi:Y\to X\times\P^1$ be a double point 
degeneration over $0\in \P^1$. Let
$$Y_0= A\cup B \to X$$
be the fiber over 0, and let 
$Y_\infty\to X$ be a regular fiber.
The associated double point relation is
\[
[Y_\infty\to X]=[A\to X]+[B\to X]-[\P(\pi)\to X]\ \in \omega_*(X).
\]

\begin{lem} \label{ddd}
 Let $L$ be a line bundle on $X$. Then,
\[
\cn(L)([Y_\infty\to X])=\cn(L)\Big([A\to X]+[B\to X]-[\P(\pi)\to X]\Big).
\]
\end{lem}

\begin{proof} 
The various classes $\cn(L)([W\to X])$ are defined by operating on 
$\omega_*(W)$ and then pushing forward to $X$. Hence,  
we may replace $X$ with $Y$, $L$ with $\pi^*p_1^*L$, and $\pi$ 
with 
$$(\id_Y,p_2\circ \pi):Y\to Y\times\P^1.$$

Since  $Y\in\Sm_k$, we may choose a very ample line bundle $M$ for
which $L\otimes M$ is also very ample. Then,
\[
\cn(L)=F^-_\omega(\cn(L\otimes M),\cn(M))
\]
is a map from $\sM_*(Y)^+$ to $\omega_{*-1}(Y)$. 
The result follows from Lemmas~\ref{lem:ChernDescent} and \ref{lem:dim}.
\end{proof}

By Lemma \ref{ddd}, the operator
$\cn(L):\sM_*(X)^+\to \omega_{*-1}(X)$ descends to
\[
\cn(L):\omega_*(X)\to \omega_{*-1}(X).
\]
Hence, we have constructed first Chern class operators on $\omega_*$
for arbitrary line bundles.

\begin{lem} \label{lem:dim2} Let $Y\in\Sm_k$, and let 
$$L_1,\ldots, L_{r>\dim_k Y}\to Y$$ be 
 line bundles. Then, $$\prod_{i=1}^r\cn(L_i)=0$$ 
as an operator on $\omega_*(Y)$.
\end{lem}

\begin{proof}
Since  $Y$ quasi-projective,
$
\cn(L_i)=F^-_\omega(\cn(L_i\otimes M),\cn(M))
$
for any choice of very ample line bundle $M$ on $Y
$ for which $L_i\otimes M$ is very ample.
Since $$F^-_\omega(u,v)=u-v\mod (u,v)^2,$$
 Lemma~\ref{lem:dim} 
implies the result.
\end{proof}

Axioms (A3), (A4), (A5) and (A8) 
for globally generated $L$ immediately imply
 these axioms for arbitrary $L$. 
Similarly, the functoriality 
of  Lemma~\ref{lem:ChernFunct} extends to arbitrary
line bundles  $L$.

\begin{prop} \label{prop:FGL3} Let $L$ and $M$ be line bundles on 
$X\in\Sch_k$. Then,
\[
\cn(L\otimes M)=F_\omega(\cn(L),\cn(M)).
\]
\end{prop}

\begin{proof} 
By the definition of Chern classes and Lemma \ref{lem:dim2},
the operator
$$F_\omega(\cn(L),\cn(M)): \omega_*(X) \to \omega_{*-1}(X)$$
is well-defined. 

Since $\omega_*(X)$ is generated by the classes $f_*(1_Y)$ for 
$$[f:Y\to X]\in \sM(X)^+,$$ 
property (A3) can be used to reduce to the case of $X\in\Sm_k$.

Take very ample line bundles $N_1$, $N_2$ on $X$
 with $L\otimes N_1$ and $M\otimes N_2$ very ample. 
Then,
 $$L\otimes M\otimes N_1\otimes N_2, \ \
N_1\otimes N_2$$ are also very ample.
 The Proposition follows from Proposition~\ref{prop:FGL2} and the power series identity
\[
F_\omega(F^-_\omega(u_1,v_1), F^-_\omega(u_2,v_2))=
F^-_\omega(F_\omega(u_1,u_2), F_\omega(v_1,v_2)),
\]
after taking
\begin{align*}
&u_1= \cn(L\otimes N_1), \ v_1=\cn(N_1),\\
&u_2=\cn(M\otimes N_2),\ v_2=\cn(N_2).
\end{align*}
\end{proof}

\subsection{Proof of Theorem \ref{two}}
Double point cobordism theory $\omega_*$ was shown in Section \ref{bmp}
to define a Borel-Moore functor with product:
structures (D1), (D2), and (D4) satisfying axioms (A1), (A2), (A6),
and (A7).

We have added first Chern classes (D3) and verified axioms
(A3), (A4), (A5), and (A8). Hence, $\omega_*$ is oriented.

The formal group law defined by Proposition \ref{formgp} yields a
canonical ring homomorphism
$$\L_* \to \omega_*(k).$$
Hence, $\omega_*$ is $\L_*$-functor.

In order for $\omega_*$ to be an oriented Borel-Moore
$\L_*$-functor of geometric type, the axioms of Section
\ref{gtype} must be satisfied. Axiom (Dim) is Lemma \ref{lem:dim2}, and
axiom (FGL) is Proposition \ref{prop:FGL3}. 
The proof of Theorem \ref{two} will be completed by
establishing the remaining axiom (Sect).

\section{Axiom (Sect)}  
\label{eee}
\subsection{The difference series} 
Since the Chern class operator $\cn(L)$ for a general line bundle 
$L$ is defined using the difference $F^-_\omega$ in our formal group law, 
we will require a universal construction of $F^-_\omega$ along the lines 
of our construction of $F_\omega$. 

The variety $Y_{n,m}$, defined in Section \ref{ind},
is the total space of the line bundle $O(1,-1)$ on $\P^n\times\P^m$ 
with projection $\pi$ and zero-section $s$,
$$\pi:Y_{n,m}\to \P^n\times\P^m, \ \ 
s:\P^n\times\P^m\to Y_{n,m}.$$
 Let $S_{n,m}\subset Y_{n,m}$ be the image of the zero section. 

For $0\leq i \leq n$ and $0,\leq j \leq m$,
a closed immersion $$Y_{i,j}\to Y_{n,m}$$ is induced by a 
choice of bi-linear embedding $\P^i\times \P^j\to\P^n\times\P^m$.

\begin{lem}\label{lem:DifferenceRelation} For $n,m\ge0$,  
\begin{equation}\label{eqn:DifEq}
[S_{n,m}\to Y_{n,m}]
=\sum_{i=0}^n \sum_{j=0}^m 
b_{i,j}^{n,m}\cdot
[Y_{n-i,m-j}\to Y_{n,m}]\ \in \omega_*(Y_{n,m})
\end{equation}
for $b^{n,m}_{i,j}\in \omega_{i+j-1}(k)$.
\end{lem}

\begin{proof} If $n=m=0$, then
$Y_{n,m}=\A^1$ with $S_{n,m}\to Y_{n,m}$ given by
the inclusion of 0. Clearly $[0\to\A^1]=0$ in $\omega_0(\A^1)$, whence the 
result.{\footnote{Consider the  morphism
 $\pi:\A^1 \to \A^1 \times \P^1$
determined by $(\id,i)$ where $$i:\A^1 \to \P^1$$ is the inclusion obtained
by omitting $0\in \P^1$.
The projective morphism $\pi$ is a double point degeneration over
$0\in \P^1$,
$$\pi^{-1}(0)= \emptyset \cup \emptyset.$$
The associated double point cobordism shows
$[\Spec(k) \to \A^1]=0$ in $\omega_*(\A^1)$ for every
closed point.}}

We proceed by induction on $(n,m)$.
We have the linear equivalence
\[
S_{n,m}+Y_{n,m-1}\sim Y_{n-1,m}
\]
on $Y_{n,m}$.
Clearly $S_{n,m}+Y_{n,m-1}+Y_{n-1,m}$ 
is a reduced strict normal crossing divisor on $Y_{n,m}$. 
By Lemma~\ref{lem:ExtendedDP}, we obtain the relation 
\begin{multline*}
[S_{n,m}\to Y_{n,m}]=[Y_{n-1,m}\to Y_{n,m}]-[Y_{n,m-1}\to Y_{n,m}]\\
+[\P_1\to Y_{n,m}]-[\P_2\to Y_{n,m}]+[\P_3\to Y_{n,m}]
\end{multline*}
where $\P_1\to S_{n,m-1}$ is an admissible $\P^1$-bundle, 
$\P_2\to S_{n-1,m-1}$ is an admissible tower, and
$\P_3\to S_{n-1,m-1}$ is an admissible $\P^2$-bundle.

We apply Lemma~\ref{lem:TowerDecomp} to  $\P_1\to S_{n,m-1}$ 
with generators $S_{n-1,m-1}$ and $S_{n,m-2}$ for $\Pic(S_{n,m-1})$.
Similarly, we apply Lemma \ref{lem:TowerDecomp}
to  $\P_2\to S_{n-1,m-1}$ and $\P_3\to S_{n-1,m-1}$. 
We find
\begin{multline*}
[S_{n,m}\to Y_{n,m}]=[Y_{n-1,m}\to Y_{n,m}]-[Y_{n,m-1}\to Y_{n,m}]\\
+\sum_{i=0}^n\sum_{j=1}^m c_{i,j}\cdot[S_{n-i,m-j}\to Y_{n,m}]
\end{multline*}
with  $c_{i,j}\in\omega_*(k)$.

Since $S_{n-i,m-j}\to Y_{n,m}$ factors 
through $S_{n-i,m-j}\to Y_{n-i,m-j}$, the induction hypothesis 
finishes the proof.
\end{proof}

For $0\leq i+j\le n$, $0\leq j\le m$, the elements $b_{i,j}^{n,m}$ on the right
side of \eqref{eqn:DifEq} are uniquely determined by 
Proposition~\ref{prop:ProductInd3}.

\begin{lem}\label{lem:DifferenceRelationUnique} 
If $N\geq n$, $M\geq m$, then
$$b_{i,j}^{n,m}=b_{i,j}^{N,M}$$
for $0\leq i+j\le n$, $0\leq j\le m$.
\end{lem}

\begin{proof} 
The bi-linear embedding $\P^n\times\P^m\to\P^N\times\P^M$ 
induces a closed embedding $$i:Y_{n,m}\to Y_{N,M}$$ 
which satisfies the conditions of the second
example of Section \ref{exx}. Thus, we have a well-defined 
$\omega_*(k)$-linear pull-back
\[
i^*:\omega_*(Y_{N,M})\to \omega_{*-d}(Y_{n,m})
\]
with $d=N-n+M-m$. Clearly 
$$i^*([S_{N,M}\to Y_{N,M}]) = [S_{n,m} \to Y_{n,m}],$$
$$i^*([Y_{N-i,M-j}\to Y_{N,M}])=[Y_{n-i,m-j}\to Y_{n,m}],$$ so 
the uniqueness statement implies the result.
\end{proof}

By Lemma \ref{lem:DifferenceRelationUnique}, we may define
$b_{i,j}\in \omega_*(k)$ by
\[
b_{i,j}=\lim_{n,m\to\infty}b^{n,m}_{i,j}.
\]
By the proof of Lemma \ref{lem:DifferenceRelation},
 $b_{0,0}=0$, $b_{1,0}=1$, and $b_{0,1}=-1$.

\begin{lem}\label{lem:DifPowerSeries} 
$F^-_\omega(u,v)=\sum_{i,j}b_{i,j}u^iv^j$.
\end{lem}

\begin{proof} 
Let $n, m\ge0$, and let $N=n+2m$, $M=m$. 
The morphism 
$$h:\P^{n+m}\times\P^m\to Y_{N,M}$$ 
was constructed in the proof of Proposition~\ref{prop:ProductInd3}. 
We see 
$$h^{-1}(S_{N,M}) = \P^{n+m-1}\times\P^{m} \to \P^{n+m}\times \P^m$$
is
a bi-linear embedding
and 
$$h^{-1}(Y_{N-i,M-j}) = H^{(i)}_{n+m,m-j}\to \P^{n+m}\times \P^m.$$
The relation \eqref{eqn:DifEq} for $(N,M)$ pulls back under $h$ to
\[
[\P^{n+m-1}\times\P^m\to \P^{n+m} \times \P^m]
=\sum_{i,j}b^{N,M}_{i,j}\cdot
[H^{(i)}_{n+m,m-j}\to \P^{n+m}\times\P^m].
\]
We have $b^{N,M}_{i,j}=b_{i,j}$ for 
$$0\leq i+j\le N=n+2m, \ \ 0\leq j\le M=m.$$ 
Since $H^{(i)}_{n+m,m-j}\to \P^{n+m}\times \P^m$ has codimension $i+j$  
and is empty if $j>m$,  
\[
[\P^{n+m-1}\times\P^m]=\sum_{i=0}^{n+m}\sum_{j=0}^mb_{i,j}\cdot
[H^{(i)}_{n+m,m-j}\to \P^{n+m}\times\P^m].
\]

Consider the formal group law determined by $\omega_*$.
The difference $F^-_\omega$ admits a power
series expansion, 
$$F^-_\omega(u,v)=\sum_{i,j}\tilde{b}_{i,j}u^iv^j,$$
where $\tilde{b}_{i,j}\in \omega_*(k)$.
Certainly,
\[
\cn(O(1,0))(1_{\P^{n+m}\times\P^m})=
F^-_\omega(\cn(O(1,1)),\cn(O(0,1))(1_{\P^{n+m}\times\P^m}).
\]
Since 
\begin{align*}
&[\P^{n+m-1}\times\P^m]=\cn(O(1,0))(1_{\P^{n+m}\times\P^m}),\\
&[H^{(i)}_{n+m,m-j}\to \P^{n+m}\times\P^m]=
\cn(O(1,1))^i\cn(O(0,1))^j(1_{\P^{n+m}\times\P^m}),
\end{align*}
we find
\[
[\P^{n+m-1}\times\P^m]=
\sum_{i,j}\tilde{b}_{i,j}\cdot [H^{(i)}_{n+m,m-j}\to \P^{n+m}\times\P^m].
\]
Therefore,
\[
\sum_{i,j}(b'_{i,j}-b_{i,j})\cdot
[H^{(i)}_{n+m,m-j}\to \P^{n+m}\times\P^m]=0.
\]
By Proposition~\ref{prop:ProductInd2},  $b_{i,j}=b'_{i,j}$ for 
$0\leq i+j\le n+m$, $0\leq j\le m$. 
As $n$ and $m$ were arbitrary, the proof is complete.
\end{proof}

\subsection{Proof of Theorem \ref{two}}
We now complete the last step in the proof of Theorem \ref{two}.

\begin{prop} \label{prop:Sect}
Double point cobordism $\omega_*$ satisfies axiom (Sect).
\end{prop}

\begin{proof} Let $Y\in \Sm_k$ be of dimension $d$.
  Let $L$ be a line bundle on $Y$ with
transverse section $s\in H^0(Y,L)$. 
Let 
$D\subset Y$ be the smooth divisor associated to $s$.

Let $M$ be a very ample line bundle on $Y$ for which
$L\otimes M$ is also very ample.
Let $$f:Y\to \P^n, \ \ g:Y\to \P^m$$
 be closed embeddings satisfying
 $$L\otimes M\cong f^*O(1),  \ \ M\cong g^*O(1).$$
Certainly, $d\le n$, $d\le m$. 

Let $h=(f,g):Y\to \P^n\times\P^m$. The section $s$ defines a 
lifting $$(h,s):Y\to Y_{n,m}$$ which satisfies the conditions of 
the second example of Section \ref{exx}.
We obtain a well-defined $\omega_*(k)$-linear pull-back
\[
(h,s)^*:\omega_*(Y_{n,m})\to \omega_{*-n-m-1+d}(Y).
\]
By construction, $(h,s)^*([S_{n,m}\to Y_{n,m}])=[D\to Y]$.

Since
$\cn(\pi^*O(1,0))^i\cn(\pi^*O(0,1)^j(1_{Y_{n,m}})=[Y_{n-i,m-j}\to Y_{n,m}]$
and 
\[
(h,s)^*(\pi^*O(1,0))=L\otimes M,\ (h,s)^*(\pi^*O(0,1))=  M,
\]
Lemma~\ref{lem:DifferenceRelationUnique},  
Lemma~\ref{lem:DifPowerSeries},  and the naturality of $\cn$ given by 
Lemma~\ref{lem:ChernFunct} yield
\[
(h,s)^*(\sum_{i,j} b_{i,j}^{n,m}[Y_{n-i,m-j}\to Y_{n,m}])=
F^-_\omega(\cn(L\otimes M), \cn(M))(1_Y).
\]
The ``error terms" arising from any inequalities  
$b_{i,j}^{n,m}\neq b_{i,j}$ vanish because
\[
(h,s)^*([Y_{n-i,m-j}\to Y_{n,m}])=0
\]
if 
$i+j>n\ge d$ or if $j>m$ for dimensional reasons. 

Applying $(h,s)^*$ to the relation \eqref{eqn:DifEq} yields the identity
\[
[D\to S]=F^-_\omega(\cn(L\otimes M), \cn(M))(1_Y) =\cn(L)(1_Y),
\]
which verifies axiom (Sect).
\end{proof}

\section{Theorem~\ref{one} and Corollary \ref{bas}} 
\label{t1}

\noindent{\em Proof of Theorem \ref{one}.}
For clarity, we write
 $[f:Y\to X]_\omega$ for 
$$[f:Y\to X]\in \omega_*(X)$$
and
$[f:Y\to X]_\Omega$ for the associated class  in $\Omega_*(X)$. 
Similarly, let 
$$1^\omega_Y=[\id_Y]_\omega,\ \ 1^\Omega_Y=[\id_Y]_\Omega.$$

By Proposition \ref{prop:Map}, 
there is natural transformation 
$$\vartheta:\omega_* \to \Omega_*$$ 
of Borel-Moore functors on $\Sch_k$,  
\[
\vartheta_X([f:Y\to X]_\omega)=[f:Y\to X]_\Omega\in\Omega_*(X). 
\]
Moreover, $\vartheta_X$ is surjective for every $X\in \Sch_k$.

By Theorems \ref{two} and \ref{thm:OBMGUniversal}, 
there is a natural transformation 
$$\tau:\Omega_*\to\omega_*$$ 
of oriented Borel-Moore functors of geometric type.
Let $Y\in \Sm_k$, and let 
$$p:Y \to \Spec(k)$$ be the
structure map.
Since
$$1_Y^\Omega = p^*(1), \ \ 1_Y^\omega= p^*(1),$$
and $\tau$ respects the unit and smooth pull-back,
$$\tau(1_Y^\Omega)= 1_Y^\omega.$$
Hence,
\begin{eqnarray*}
\tau_X([f:Y\to X]_\Omega)&=&\tau_X(f_*(1_Y^\Omega))\\
&=& f_*(\tau_Y(1_Y^\Omega))\\
& = & f_*(1_Y^\omega) \\
&= &[f:Y\to X]_\omega.
\end{eqnarray*}
Therefore $\tau\circ\vartheta=\id_\omega$, so $\vartheta$ is an isomorphism.
\qed

\vspace{10pt}

\noindent{\em Proof of Corollary \ref{bas}.} 
We may assume $k\subset \mathbb{C}$.
The canonical homomorphism $$\Omega^*(k) \to MU^{2*}(\text{pt})$$
discussed in Section \ref{qui} is an isomorphism.
Since $MU^{2*}(\text{pt})$ is well-known to have a
rational basis determined by the products of projective
spaces, the Corollary is deduced from Theorem \ref{one}.
\qed

\section{Donaldson-Thomas theory}
\label{last}

\subsection{Proof of Conjecture 1}
Let $\mathbb{Q}[[q]]^*\subset \mathbb{Q}[[q]]$ 
denote the multiplicative group of power
series with constant term 1.
Define a group homomorphism 
$$\mathsf{Z}: (\sM_3(\Spec(\C))^+,+) \to (\mathbb{Q}[[q]]^*,\cdot)$$
on generators
by the partition function for degree 0 Donaldson-Thomas theory 
defined in Section \ref{dt},
$$\mathsf{Z}([Y]) = \mathsf{Z}(Y,q).$$
We use here the abbreviated notation 
$$[Y]=[Y\to \Spec(\C)]\in \sM_3(\Spec(\C)).$$
Since double point relations  
hold in Donaldson-Thomas
theory \eqref{dpdt}, the homomorphism
$\mathsf{Z}$ descends to $\omega_*(\C)$,
$$\mathsf{Z}: \omega_*(\C) \to \mathbb{\Q}[[q]]^*.$$

By Corollary \ref{bas}, the class $[Y]\in \omega_3(\C)$
is expressible rationally in terms of the classes
$$[\P^3], \ [\P^2 \times \P^1], \ [\P^1\times \P^1\times \P^1].$$
Hence, 
$$r[Y] = s_3 [\P^3] + s_{21} [\P^2\times \P^1] + s_{111} 
[\P^1\times\P^1\times\P^1] \ \in \omega_*(\C)$$
for integers $r\neq 0$, $s_3$, $s_{21}$, and $s_{111}$.
Therefore
\begin{equation}
\label{cvf}
\mathsf{Z}(Y,q)^r = 
\prod_{|\lambda|=3} \mathsf{Z}(\P^\lambda,q)^{s_\lambda}.
\end{equation}

Conjecture 1 has been proven for 3-dimensional products
of projective spaces in \cite{MNOP1,MNOP2}. The right side of
\eqref{cvf}
can therefore be evaluated:
\begin{eqnarray*}
\prod_{|\lambda|=3} \mathsf{Z}(\P^\lambda,q)^{s_\lambda} & = &
\prod_{|\lambda|=3} M(-q)^{s_\lambda \int_{\P^\lambda} c_3(T_{\P^\lambda}
\otimes K_{\P^\lambda})} \\
& = &
M(-q)^{\sum_{|\lambda|=3} 
s_\lambda \int_{\P^\lambda} c_3(T_{\P^\lambda}
\otimes K_{\P^\lambda})}
\end{eqnarray*}
Since algebraic cobordism respects Chern numbers \footnote{For example,
because complex cobordism does.},
$$\mathsf{Z}(Y,q)^r =
M(-q)^{ 
r \int_{Y} c_3(T_{Y}
\otimes K_{Y})}.$$
Finally, since $Z(Y,0)=1$ and $M(0)=1$,
$$\mathsf{Z}(Y,q) =
M(-q)^{ 
\int_{Y} c_3(T_{Y}
\otimes K_{Y})},$$
completing the proof.\qed

\subsection{Conjecture $1'$}
Next, we consider an equivariant version of Conjecture 1 proposed in
\cite{bp}.

Let $X$ be a smooth quasi-projective 3-fold over $\C$ equipped
with an action of an algebraic torus $T$ with {\em compact}
fixed locus $X^T$. If $X^T$ is compact, 
$\text{Hilb}(X,n)^T$ is also compact, and
$$N^X_{n,0} = \int_{[\text{Hilb}(X,n)^T]^{vir}} 
\frac{1}{e(\text{Norm}^{vir})} \ \in \Q(\mathbf{t})$$
is well-defined \cite{GP}. 
Here 
 $$\mathbf{t}=\{t_1, \ldots, t_{\text{rk}(T)}\}$$
is a set of  generators
of the $T$-equivariant 
cohomology of a point. 
Let 
$$\ZZ(X,q,\mathbf{t}) = 1 + \sum_{n\geq 1} N^X_{n,0}\ q^n$$
be the equivariant partition function.

Since $X^T$ is compact, the right side of the equality of 
Conjecture 1 is also well-defined via localization,
$$\int_X c_3(T_X \otimes K_X) = 
\int_{X^T} \frac{c_3(T_X\otimes K_X)}{e(\text{Norm})} 
\in \Q(\mathbf{t}).$$

\vspace{10pt}

\noindent{\bf Conjecture $1'$.} \cite{bp}\,  $\ZZ(X,q,\mathbf{t}) =
M(-q)^{\int_X c_3(T_X \otimes K_X)}.$

\vspace{10pt}

We will prove Conjecture $1'$ before proving Conjecture 2 for 
relative Donldson-Thomas theory.

\subsection{Local geometries}
Let $M$  be a smooth projective variety over $\C$
of pure dimension at most 3.
Let 
$$N\to M$$
be a vector bundle of satisfying
$$\text{rk}(N)= 3-\dim_\C M.$$

The space total space $N$
may be viewed as a {\em local} neighborhood\footnote{There
is no algebraic tubular neighborhood result even formally.}  of $M$
in a 3-fold embedding. 
If $$N =\bigoplus_{i=1}^r N_i$$
is a direct sum decomposition, an $r$-dimensional 
torus $T$ acts canonically on the total space $N$ 
by scaling the factors of $N$.
Since $N^T=M$, the fixed locus is compact.

We will first prove Conjecture $1'$ for the local geometry $N$.
In case $M$ has dimension $0$ or $1$, Conjecture $1'$ has
been proven in \cite{MNOP1,MNOP2} and \cite{OP} respectively.
If $Y$ has dimension $3$, Conjecture $1'$ reduces
to Conjecture $1$. Only the dimension 2 case remains.

\subsection{Proof of Conjecture $1'$ for local surfaces}
The proof relies upon 
a double point cobordism theory for local geometries.
To abbreviate the discussion, we focus our attention on
the double point cobordism theory for local surfaces
over $\Spec(\C)$.

Consider the free group $M_{2,1}(\C)^+$ generated
by pairs $[S,L]$ where $S$ is smooth, irreducible, projective surface and
$$L\to S$$ is
a line bundle.
The subscript $(2,1)$ captures the dimension of $S$ and the
rank of $L$.
We define a double point cobordism theory
$\omega_{2,1}(\C)$ as a quotient of $M_{2,1}(\C)^+$
by double point relations.

Double point relations are easily defined in the local setting.
Let 
$$\pi: \mathcal{S} \to \P^1$$
be a projective morphism determining  
a double point degeneration with
$$\mathcal{S}_0= A\cup B,$$ and let 
$$\mathcal{L} \to \mathcal{S}$$
be a line bundle.
For each regular value $\zeta\in\P^1$ of $\pi$, define an
associated double point relation by
\begin{equation}\label{ddpptt}
[\mathcal{S}_\zeta,\mathcal{L}_\zeta] -[A,\mathcal{L}_A]-[B,
\mathcal{L}_B]+ [\P(\pi),\mathcal{L}_{\P(\pi)}].
\end{equation}
Here, subscripts denote restriction (or, in the case of 
$\mathcal{L}_{\P(\pi)}$, pull-back).

Let $\mathcal{R}_{2,1}(\C) \subset M_{2,1}(\C)^+$
be the subgroup generated by all double point relations. 
Double point cobordism theory for local surfaces
is defined by
$$\mathcal{\omega}_{2,1}(\C) = M_{2,1}(\C)^+/\mathcal{R}_{2,1}(\C).$$

\begin{lem} Double point cobordism theory \label{genn} 
$\omega_{2,1}(\C)$ for local surfaces
is generated (over $\Q$) by elements of the following form:
\begin{enumerate}
\item[(i)] $[\P^2, O_{\P^2}]$,
\item[(ii)] $[\P^1\times \P^1, L]$,
\item[(iii)] $[F_1, L]$,
\end{enumerate}
where $F_1$ is the blow-up of $\P^2$ in a point. 
\end{lem}

\begin{proof}
There is a natural group homomorphism
$$\iota:\omega_2(\C)\otimes_\Z \Q \to \omega_{2,1}(\C)\otimes_\Z \Q$$
defined by $\iota([S]) = [S,O_S]$. By Corollary \ref{bas},
the image of $\iota$ is generated by
$$[\P^2,O_{\P^2}], \ \ [\P^1\times \P^1, O_{\P^1\times \P^1}].$$

Let $[S,O_S(C)]\in M_{2,1}(\C)^+$ where $C\subset S$
is smooth divisor.
Consider the deformation to the normal cone of $C\subset S$,
$$\pi:\mathcal{S} \to \P^1$$
with degenerate fiber
$$\mathcal{S}_0= S \cup \P(O_C\oplus O_C(C)).$$
Since $\mathcal{S}$ is the blow-up of $S\times \P^1$ along
$C \times 0$, there is a canonical morphism
$$\nu: \mathcal{S} \to S$$
obtained from blow-down and projection.
Let $\mathcal{L} \to \mathcal{S}$
be defined by
$$\mathcal{L} = \nu^*(O_S(C+D)) \otimes O_{\mathcal{S}}(-\P(O_C\oplus
O_C(C))).$$
where $D$ is a Cartier divisor on $S$. 
The double point
relation associated to $\mathcal{L}\to \mathcal{S}$ is
\begin{equation}\label{vtr}
[S,O_S(C+D)] - [S,O_S(D)]-[\P(O_C\oplus O_C(C))
,L'] + [\P(\pi),L'']
\end{equation}
where $L'$ and $L''$ are line bundles.

Let $\Gamma\subset \omega_{2,1}(\C)$ be the subgroup generated
by $\text{Im}(\iota)$
and elements of the form
$
[P,L]
$
where $P$ is a $\P^1$-bundle over a smooth projective
curve.
If $D$ is taken to be 0 in \eqref{vtr}, we find
$[S,O_S(C)] \in \Gamma$.
For general a Cartier divisor $D$, 
$$[S,O_S(C+D)]\in \Gamma \iff [S,O_S(D)]\in \Gamma.$$
 Since, for any $D$, there
exists smooth curves $C,C'$ for which 
$$O_S(C+D)\cong O_S(C'),$$
we find
$\Gamma=\omega_{2,1}(\C)$.

By elementary degenerations, elements of type (ii) and (iii)
generate the classes $[\P,L]$ of $\Gamma$.
\end{proof}

The computation of the degree 0 equivariant vertex in
\cite{MNOP1,MNOP2} proves Conjecture $1'$ for the toric
generators (i-iii) of Lemma \ref{genn}.
Conjecture $1'$ then follows for local surfaces by an argument
parallel to the proof of Conjecture 1. \qed

\subsection{Proof of Conjecture $1'$}
Let $T$ be an $r$-dimensional torus acting
on a smooth quasi-projective 3-fold $X$
with compact fixed locus $X^T$. 
The 1-dimensional subtori of $T$ are described by elements of the
lattice
$\mathbb{Z}^r$. Since 1-dimensional tori
$T_1\subset T$
with equal fixed loci $$X^{T_1} = X^T$$
determine a Zariski dense subset of $\mathbb{Z}^r$, Conjecture $1'$
is implied by the rank 1 case.

We assume $T$ is a 1-dimensional torus. If the $T$-action
on $X$ is trivial, Conjecture $1'$ reduces to Conjecture 1.
We assume the $T$-action is nontrivial.
The components of the fixed locus
$$X^T = \bigcup_i X_i^T$$
are of dimension 0, 1, or 2. 
Certainly
\begin{equation}\label{gty}
\mathsf{Z}(X,q,t) = \prod_i \mathsf{Z}(X_i,q,t)
\end{equation}
where 
$$\mathsf{Z}(X_i,q,t) = 
\sum_n q^n \int_{[\text{Hilb}(X,n)_i^T]^{vir}} 
\frac{1}{e(\text{Norm}^{vir})}$$
and $\text{Hilb}(X,n)_i^T\subset \text{Hilb}(X,n)^T$
is locus supported on $X_i^T$.
We will prove
\begin{equation}\label{gtyy}
\mathsf{Z}(X_i,q,t)= M(-q)^{ \int_{X_i^T}  
\frac{c_3(T_X\otimes K_X)}{e(\text{Norm}_i)}}
\end{equation}
where $\text{Norm}_i$ is the normal bundle of
$X_i^T \subset X$. Conjecture $1'$
follows from \eqref{gty} and \eqref{gtyy}.

Equality \eqref{gtyy} is proven separately
for each possible dimension of  $X_i^T$. The dimension 1
case is the most delicate.

\vspace{10pt}
\noindent{\bf Dim 0}. If $X_i^T=p$ is a point, then by Theorem
2.4 of \cite{BB}, the $T$-action on $X$ is analytically equivalent
in a Euclidean neighborhood of $p$
to the $T$-action on the tangent space $T_p(X)$.  
The $T$-action at a point $u\in U$ of a
Euclidean neighborhood is defined only
locally at $1\in T$.
Equality
\eqref{gtyy} in the dimension 0 case
follows from
the degree 0 vertex evaluation of \cite{MNOP1,MNOP2}.

\vspace{+10pt}
\noindent{\bf Dim 2}. If $X_i^T=S$ is a surface, 
the $T$-weight on the normal bundle of $S\subset X$
may be assumed positive. 
The Bialynicki-Birula stratification \cite{BB} provides
a $T$-equivariant
Zariski neighborhood
of $S$ determined by a $T$-equivariant
affine bundle 
$$S_+ \to S$$
of rank 1 with a $T$-fixed section.
In the rank 1 case,
$S_+$ is the total space of a $T$-equivariant
line
bundle over $S$. 
Equality
\eqref{gtyy} in the dimension 2 case
follows from
Conjecture $1'$ for local surfaces.

\vspace{+10pt}
If $X_i^T=C$ is a curve, there are three possibilities.
Let $N_C$ be the rank 2 normal bundle of $C\subset X$.
The $T$-representation on the fiber of $N_C$ has
nontrivial weights $w_1$ and $w_2$. 

\vspace{+10pt}
\noindent{\bf Dim 1, weights of opposite sign}. If the weights $w_1$ and
$w_2$ have { opposite} signs, then there is a canonical
$T$-equivariant splitting 
$$N=N_+\oplus N_-$$
as a sum of line bundles.
The Bialynicki-Birula
stratification yields quasi-projective surfaces 
$$C_+,C_-\subset X$$ 
corresponding to the positive and negative normal directions.
Since the affine bundles 
$$C_\pm\to C$$
are of rank 1 with $T$-fixed sections, there
 are $T$-equivariant isomorphisms
$$\phi_\pm : C_\pm \to N_\pm$$ 
where the total spaces of the line bundles occur on the right.

Let $p\in C$.
By Theorem 2.4 of \cite{BB}, the $T$-action on a Euclidean
neighborhood $U_X\subset X$ of $p\in X$ is analytically equivalent to the
$T$-action on a Euclidean neighborhood $U_N \subset N_C$ 
of $p\in N_C$. 
Certainly the images of $C_\pm$ are the intersections of
$U$ with $N_\pm$.

Since the $T$-action on $N_C$ has weights of opposite sign,
the $T$-equivar\-iant automorphism group of $U$ over $C$ 
which fixes $U\cap N_\pm$ pointwise is trivial.
In particular, there is a unique $T$-equivariant isomorphism
$$U_X \to U_N$$
compatible with $\phi_\pm$. 
Patching together the isomorphisms yields an $T$-equivariant
analytic
isomorphism between $X$ and $N_C$ defined in a Euclidean
neighborhood of $C$. 
Equality
\eqref{gtyy} in the 1-dimensional opposite sign case then 
follows from
Conjecture $1'$ for local curves proved in \cite{OP}.

\vspace{+10pt}
If the weights $w_1$ and $w_2$ are of the same sign, we may
assume the weights to be positive.
The Biaylnicki-Birula stratification yields a $T$-equivariant
Zariski neighborhood
of $C$ determined by a $T$-equivariant
affine bundle
$$C_+\to C$$
of rank 2. We will see $C_+$ need not be the
total space of a $T$-equivariant rank 2 vector bundle on $C$.

The weights $w_1$ and $w_1$ are {\em related} if there
exists an integer $k\geq 2$ for which either
$$w_1\cong kw_2 \ \ \text{or}\ \ kw_1 \cong w_2.$$

\vspace{+10pt}
\noindent{\bf Dim 1, related weights of same sign}.
Without loss of generality,
we may assume the relation is $w_1=kw_2$.

Let $\mathbb{C}^2$ be a $T$-representation with weights 
$w_1$ and $w_2$,
$$t\cdot(z_1,z_2) = (t^{w_1}z_1,t^{w_2} z_2).$$
The $T$-equivariant
automorphism group $G$ of $\mathbb{C}^2$ is given by 
$2\times 2$ upper triangular matrices,
\begin{equation}
\label{ftyr}
\gamma_{\left( \begin{array}{cc}
 \lambda_1&    \delta  \\
    0 &  \lambda_2 \end{array}\right)}\Big(z_1,z_2\Big)
=\Big(\lambda_1 z_1 + \delta z_2^k,\lambda_2 z_2\Big).
\end{equation}

Every Zariski locally trivial $G$-torsor $\tau$ on $C$
yields an $T$-equivariant affine bundle 
$$A_\tau \to C$$
of rank 2 over
$C$ with a $T$-equivariant section.
The bundle $A_\tau$ is obtained by the $G$-action \eqref{ftyr}.
The family of homomorphisms
$$\rho_\xi: G \to G$$
for $\xi \in \mathbb{C}$ defined by
$$\rho_\xi \left( \begin{array}{cc}
 \lambda_1&    \delta  \\
    0 &  \lambda_2 \end{array}\right) =
\left( \begin{array}{cc}
 \lambda_1&    \xi\cdot \delta  \\
    0 &  \lambda_2 \end{array}\right)$$
is a algebraic deformation of the identity $\rho_1$ to the
the diagonal projection
$$\rho_0: G \rightarrow (\mathbb{C}^*)^2.$$
For each $G$-torsor $\tau$, let $\tau_\xi$ be the $G$-torsor
induced by $\rho_\xi$. Then, the algebraic family $A_{\tau_\xi}$
of $G$-torsors is a $T$-equivariant deformation of $A_\tau$
to $A_{\tau_0}$. The latter is the total space of
a $T$-equivariant vector bundle on $C$.

Bialynicki-Birula proves the $T$-equivariant
affine bundle
$$C_+\to C$$
is obtained from a $G$-torsor as above. Since $C_+$ is
$T$-equivariantly deformation equivalent to the total
space of a rank 2 vector bundle over $C$, equality
\eqref{gtyy} follows from the local curve case together with
the deformation invariance of the virtual class.

\vspace{+10pt}
\noindent{\bf Dim 1, unrelated weights of the same sign}.
If $w_1$ and $w_2$ are not related, 
$$C_+ \to C$$
is the total space of a $T$-equivariant rank 2 vector
bundle over $C$, see Section 3 of \cite{BB}.
Equality
\eqref{gtyy} then
follows from
Conjecture $1'$ for local curves. \qed

\vspace{+10pt}

\subsection{Proof of Conjecture 2}
Let $X$ be a smooth projective 3-fold over $\C$, and let
Let $S\subset X$ be a smooth surface.
Let 
$$\P=\P(O_S \oplus O_S(S)).$$
Let $S_+,S_-\subset \P$ denote
 the sections with respective normal bundles
$O_S(S)$, $O_S(-S)$ corresponding to the quotients 
$O_S(S)$, $O_S$.

We will study the Donaldson-Thomas theory of $\mathsf{\P}/S_-$
by localization.
A 1-dimensional scaling torus $T$ acts on $\P$ with  
$$\P^T= S_+ \cup S_-$$
and normal weights $t$ and $-t$ along $S_+$ and $S_-$
respectively. 
The components of the
$T$-fixed loci of $I_n(\P/S_-,0)$ lie over either $S_-$ or $S_+$.

A Donaldson-Thomas theory of {\em rubber} naturally arises on the 
fixed loci of $I_n(\P/S_-,0)$ over $S_-$.
Let
$$\mathsf{W}_- = 1+\sum_{n\geq 1} q^n
\int_{[I_n(\P/S_-\cup S_+, 0)\ \tilde{}\ ]^{vir}}
\frac{1}{-t-\Psi_+}$$
denote the rubber contributions. Here, $I_n(\P/S_-\cup S_+,0)\ \tilde{}$
denotes the rubber moduli space, and $\Psi_+$ denotes the
cotangent line associated to target degeneration.
However, since the virtual dimension of the rubber
space $I_n(\P/S_-\cup S_+,0)\ \tilde{}$
is $-1$,
$$\mathsf{W}_- = 1.$$
A discussion of virtual localization in
relative Donaldson-Thomas theory and rubber
moduli spaces can be found in \cite{MNOP2}.
See \cite{mp} for a construction of $\Psi_+$.

A local neighborhood of
$S_+\subset \P$ is given by the total space 
$$\P_+=\P\setminus S_-$$ of the line bundle
$$O_S(S) \to S_+.$$
Hence, the contributions over $S_+$ are determined
by Conjecture $1'$ for local surfaces,
$$\mathsf{W}_+ = M(-q)^{\int_{\P_+} c_3(T_{\P_+} \otimes K_{\P_+})}.$$
The equivariant integral in the exponent is easily computed
$$
{\int_{\P_+} c_3(T_{\P_+} \otimes K_{\P_+})} =
\int_\P c_3(T_\P[-S_-] \otimes K_\P[S_-]).$$

The product of the localization contributions over $S_-$ and $S_+$
yields the partition function,
\begin{eqnarray*}
\mathsf{Z}(\P/S_-,q) & = & \mathsf{W}_-  \cdot \mathsf{W}_+ \\
                     &=  &  
M(-q)^{\int_{\P} c_3(T_{\P}[-S_-] \otimes K_{\P}[S_-])}.
\end{eqnarray*}
Conjecture 2 for $\P/S_-$ is proven.

Deformation to the normal cone of $S\subset X$ yields
\begin{equation}\label{dnc}
\mathsf{Z}(X/S,q)=\mathsf{Z}(X,q)\cdot \mathsf{Z}(\P/S,q)^{-1}.
\end{equation}
Then, Conjecture 1 for $\mathsf{Z}(X,q)$ and Conjecture 2 for
$\mathsf{Z}(\P/S,q)$ imply Conjecture 2 for $\mathsf{Z}(X/S,q)$.
\qed

\vspace{+10 pt}
\noindent
Department of Mathematics \\
Northeastern University \\
360 Huntington Ave.,\\
Boston, MA 02115, USA\\
marc@neu.edu \\

\vspace{+10 pt}
\noindent
Department of Mathematics\\
Princeton University\\
Princeton, NJ 08544, USA\\
rahulp@math.princeton.edu

\end{document}